\documentclass[article]{siamart190516}%
\usepackage{bm}
\usepackage{eurosym}
\usepackage{amsmath,nccmath}
\usepackage{amsfonts}
\usepackage{amssymb}
\usepackage{graphicx}
\usepackage{fullpage}
\usepackage{url}
\usepackage{caption}
\usepackage{subcaption}
\usepackage{cite}
\usepackage{cleveref}

\newsiamremark{remark}{Remark}
\newsiamthm{conjecture}{Conjecture}
\newsiamthm{reducedmodel}{Reduced Model}
\newsiamthm{Experiment}{Experiment}
\newcommand{\br}{\mathbf{r}}

\title{Dynamics of N-spot rings with oscillatory tails in a three-component reaction-diffusion system}
\author{Yasumasa Nishiura\thanks{Research Center of Mathematics for Social Creativity, Research Institute for Electronic Science,
Hokkaido University, Sapporo, 060-0812, and Chubu University Academy of Emerging Sciences, Chubu University , Kasugai, 487-8501, Aichi, Japan.   (\email{yasumasa@pp.iij4u.or.jp}).}  \and Shuangquan Xie\thanks{School of Mathematics, Hunan University, Changsha, 410082, China (\email{xieshuangquan2013@gmail.com}).} }

\begin{document}

\maketitle

\begin{abstract}
In two-dimensional space, we investigate the slow dynamics of multiple localized spots with oscillatory tails in a specific three-component reaction-diffusion system, whose key feature is that the spots attract or repel each other alternatively according to their mutual distances, leading to rather complex patterns. One fundamental pattern is the ring pattern, consisting of $N$ equally distributed spots on a circle with a certain radius. Depending on the parameters of the system, stationary or moving (i.e., traveling and rotating)  $N$-spot rings can be observed. In order to understand the emergence of these patterns, we describe the dynamics of $N$ spots by a set of reduced ordinary differential equations (ODEs) encoding the information of each spot's location and velocity. On the basis of the reduced system, we analytically study the existence and stability of stationary and moving N-spot ring solutions, which keep most of the essential features of the collective motion of self-propelled particles. Numerical simulations of both partial differential equations (PDEs) and ODEs are provided to verify our results and the comparison between them implies several future challenges including emergence of new spots and effects of higher order terms.
\end{abstract}
\section{Introduction} 
Localized patterns are ubiquitous in nature, such as stripes in animal skin \cite{meinhardt1982models}, nerve pulses in biological systems \cite{nagumo1962active,fitzhugh1961impulses}, concentration drops of chemical reagents in chemical systems \cite{vanag2004stationary}, intensity bulbs in optical systems \cite{schapers2000interaction}, and current filaments in gas-discharge systems \cite{scholl2001nonlinear,brauer2000traveling}. The modelling of these phenomena often generates nonlinear reaction-diffusion equations that admit spatially inhomogeneous solutions localized in small regions. Investigations of many experimental and theoretical systems have shown that these localized structures' behaviors are similar to particles, exhibiting phenomena like scattering, the formation of molecule-like bound states, sef-replication, or annihilation.

 Spots are a representative class of localized structures arising in two-dimensional (2D) reaction-diffusion (RD) systems, similar to pulses or spikes in 1D. Stationary spot solutions have been reported and studied in many RD systems with suitable values of the parameters \cite{van2011planar,wei2001pattern,lloyd2009localized,bouzat1998nonequilibrium}.  For RD systems characterized by an exponentially weak spot-spot interaction, whether a stable multi-spot pattern can be observed hinges on the far field behavior of a single spot solution. The interaction of neighbouring spots is typically repulsive if the spot has monotone tails \cite{ei1994equation,ei2002pulse,ohta2001pulse}. Spots slowly drift apart, and stationary multi-spot patterns are unstable in the absence of boundaries or appropriate inhomogeneities of the control parameters. For spots with oscillatory tails, on the other hand, the spot-spot interaction oscillates with their mutual distance, allowing for infinitely many bounded states consisting of an arbitrary number of spots \cite{schenk1998interaction}. These bounded states may become unstable due to the drift bifurcation, leading to the transition from stationary spot solutions to traveling spot solutions \cite{or1998spot,gurevich2004drift}.  In general, there are no rigorous result for the existence and stability of the localized patterns such as pulse and spot with oscillatory tails. There is one exception for the FitzHugh-Nagumo equation in \cite{carter2015fast}, but this is not for three-component system and only 1D case only.  Moreover it is different from our pulse, i.e., one-handed oscillation. Stability analysis of this pulse was done in \cite{carter2016stability}. One interesting non-stationary pattern is the $N$-spot ring pattern that can be observed when $N$ spots moving slowly toward one point. With suitable initial condition, an off-center collision between spots can give rise to a stable rotating $N$-spot ring. It should be noted that a rotating two-spot ring is one of the generic patterns when started from a general initial condition of many traveling spots in a bounded domain via collision. This observation motivates us to consider the dynamics of the $N$-spot rings. Although the slow dynamics of multiple spots near the drift bifurcation can be described by a set of ODE systems \cite{zelik2009multi,ei2006interacting}, the motion of a stable multi-spot solution has not been further analyzed except in some numerical investigations \cite{liehr2003transition,liehr2004rotating}.  We aim to analytically investigate the stationary and moving $N$-spot ring patterns, as well as their stability, through the reduced ODE system.

A typical system to study an $N$-spot bound state is the following three-component RD system in $ \mathbb{R}^2$ originally introduced in \cite{schenk1997interacting} to qualitatively describe d.c. gas discharge experiments, 
\begin{equation}\label{3-rd}
\begin{aligned}
    u_t&= D_u \Delta u + k_1 u - u^3 - k_3 v - k_4 w + \kappa, \\
    \tau v_t &= D_v \Delta v + u -v,\\
    \theta w_t & = D_w \Delta w + u - w
\end{aligned} \quad \br\in  \mathbb{R}^2,
\end{equation}
where $\Delta$ is the Laplacian; $u,v,w$ depend on time $t$ and space $\br=[x,y] \in  \mathbb{R}^2 $; $k_1,~k_3,~k_4$ and $\kappa$ are kinetic parameters; the reaction rate $\tau,~\theta$ and the diffusion coefficients  $D_u,~D_v,~D_w$ are positive constants. As we are concerned with the dynamics of spots, we select a minimal model exhibiting the moving spot solutions with oscillatory tails. We consider the system \cref{3-rd} in the singular limit $D_v=0$ and $\theta=0$, equivalent to the following two-component RD system with nonlocal coupling term
\begin{equation}\label{rd}
\begin{aligned}
    u_t&= D_u \Delta u + k_1 u - u^3 - k_3 v - k_4 \mathcal{G} ^{-1} u+ \kappa, \\
    \tau v_t &=  u -v,
\end{aligned} \quad \br\in  \mathbb{R}^2,
\end{equation}
where  $\mathcal{G}$ is the operator
\begin{equation}
   \mathcal{G} :=  -D_w \Delta  + 1 
\end{equation}
For the above system, it is possible to find solutions in the form of spots with oscillatory tails, see \cite{schenk1998interaction}.  Taking $\tau$ as a bifurcation parameter, it has been shown in \cite{bode2002interaction} that $\tau_c:=1/k_3$ is the point of bifurcation from stable stationary spots to travelling ones. Above the threshold, a uniform rotating $2$-spot cluster and a traveling cluster are reported in \cite{liehr2003transition}.  It is worth noting that the order of translational and rotational bifurcations of the original PDEs for a two-spot bound state depends on the parameters $\theta$ and $D_v$, see \cite{liehr2004rotating}. Under the parameter setting $D_v=0$ and $\theta=0$, the onsets of translational and rotational bifurcations occur at the same $\tau_c$. Thus, the dynamics of the bound state beyond the bifurcation point is determined by the translational and rotational modes as well as by their interaction.  Analytically, when $\tau$ is near $\tau_c$, the set of partial differential equations can be reduced to a simple set of ordinary differential equations describing the dynamics of spots in terms of center coordinates and amplitudes of certain propagator modes.  In \cite{bode2002interaction,ei2006interacting}, with the assumption that the minimal distance between two spots is large enough, the authors derived a simplified model describing the evolution of $N$ spots near the drift bifurcation.  We summarize the reduced equations as follows:

\begin{figure}[!htb]
         \centering
         \includegraphics[width=\textwidth]{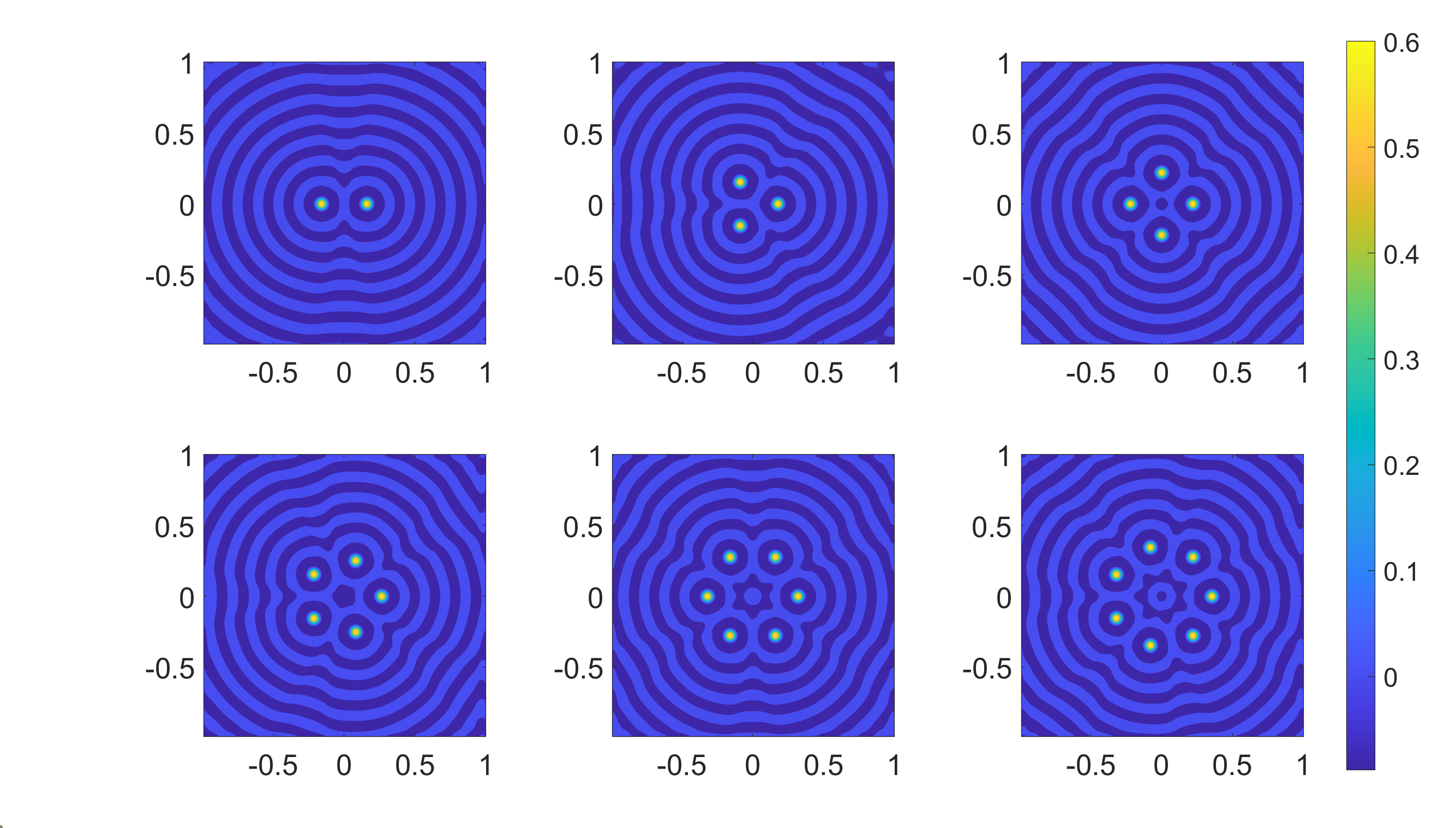}
         \caption{ Contours $u-u_c$ of stationary N-spot rings computed at $\tau=0.1$. \cref{rd} are solved with parameter $D_u=1.1\times 10^{-4},~D_w=9.64\times 10^{-4},~k_1=1.01,~k_3=0.3,~k_4=1,~\kappa=-0.1$ on a domain $[-1,1]\times[-1,1]$ with periodic boundary conditions. The initial spots are placed on a ring with the second binding radius according to the root of \cref{r0}. Stationary $N$-spot rings with the second binding radius are always stable when the parameter $\tau=0.1$, in consistent with \cref{prop1}. }
         \label{N-spot}
\end{figure}

\begin{reducedmodel}\label{lem:1}
For an N-spot ensemble in a homogeneous medium, let $p_k=x_k+i\, y_k$, where $(x_k,y_k)$ is the center of the $k$-th spot, the leading order dynamics of $N$ spots when $\tau < \tau_c$ can be described by the following ODE system:
\begin{equation}\label{eqhh0}
    \dot{p}_k=- \frac{1}{1-\tau k_3} \sum_{j\neq k} (p_k-p_j) f(|p_k-p_j|) ,
\end{equation}
where $f$ represents the interaction between spots with the form of 
\begin{equation}\label{f}
    f(d)=M_0\frac{ e^{-\alpha d}}{d^{\frac{3}{2}}} \cos(\beta (d-d_0)).
\end{equation}
with constants $M_0,~\alpha,~\beta,~d_0$ to be determined by the profile of a single spot.
\end{reducedmodel}

\begin{remark}
    The distance between two neighboring spot has to be larger than some $d_b$ that can be seen as the core of the spot, otherwise either coalescence or annihilation occurs. One necessary condition for system \cref{eqhh0} to be valid is $ |f(d)| \ll 1 $ for $d>d_b$. As the system \cref{eqhh0} diverge when $\tau=\tau_c$, we also require that $|\tau-\tau_c| \sim  \mathcal{O}(1) \gg f(d)  $.  
\end{remark}

\begin{reducedmodel}\label{lem:2}
For an N-spot ensemble in a homogeneous medium, let $p_k=x_k+i\, y_k$ and $q_k=\xi_k+i\, \eta_k$, where $(x_k,y_k)$ is the center of the $k$-th spot and $(\xi_k,\eta_k)$ is the amplitude of associated propagator mode, the leading order dynamics of $N$ spots near the drift bifurcation point $(|\tau -\tau_c| \sim f(d) \ll 1)$ can be described by the following ODE system
\begin{subequations}\label{eqhh}
\begin{equation}\label{eqhh1}
    \dot{p}_k=q_k- \sum_{j\neq k} (p_k-p_j) f(|p_k-p_j|) ,
\end{equation}
\begin{equation}\label{eqhh2}
    \dot{q}_k=M_1 q_k - M_2 q_k |q_k|^2 -k_3 \sum_{j\neq k} (p_k-p_j) f(|p_k-p_j|) ,
\end{equation}
\end{subequations}
where $M_1:=k_3^2(\tau-\frac{1}{k_3}),~M_2:= \frac{Q}{k_3} $ with constant $Q$ to be determined by the profile of a single spot and $f$ is defined in \cref{f}.
\end{reducedmodel}

It is worth noting that only the profile of a single spot is needed to evaluate constant $Q$ and function $f$ for arbitrary distances $d$, and a single spot can be computed by solving a one dimensional system since we are studying objects with locally radial symmetry. Although these simple models can describe the slow dynamics of $N$ spots, as far as the authors know, no further analytic results have been developed to study the potential stable $N$-spot configurations except in some numerical investigations \cite{liehr2003transition,liehr2004rotating,Furukawa}. On the other hand, the system \cref{eqhh} resembles the D'Orsogna model that simulates the collective motion of animal flocks  \cite{d2006self}. Several asymptotic behaviors for the D'Orsogna system arise in two dimensions, including flocking patterns and milling patterns consisting of particles distributed on a ring \cite{bertozzi2015ring,albi2014stability,kolokolnikov2011stability}. Motivated by flock and mill ring solutions emerging in the D'Orsogna model, we analyze such kinds of solutions in the reduced ODE system \cref{eqhh} and study the existence and stability of traveling and rotating ring solutions for the system \cref{eqhh}. 

Our main results are as follows:

\begin{proposition}\label{prop1}
There exist an $N$-spot ring solution to the system \cref{eqhh0} with $p_k=r_0e^{i\frac{2k\pi}{N}}$. The radius $r_0$ of the ring  is given by the root of
\begin{equation}\label{ssr}
    F(r_0)=0
\end{equation}
where
\begin{equation}
    F(r_0):=\sum_{l=1}^{N-1} (1-e^{i\theta_l })  f( 2r_0|\sin{\frac{\theta_l}{2}}|),\quad \theta_l:=\frac{2l\pi}{N}.
\end{equation}
The stability of such an $N$-spot ring solution is determined by the eigenvalues of the following $2 \times 2$ matrix for all $m=1,\cdots,N$,
 \begin{equation}
     G(m)= 
   \begin{pmatrix}
 -I_1(m)  &  -I_2(m) \\
   -I_2(m) & -I_1(-m) 
    \end{pmatrix},
 \end{equation}
where $I_1(m)$ and $I_2(m)$ are defined as
\begin{equation}
    I_1(m):=\sum_{l=1}^{N-1} \left[ f(2r_0|\sin{\frac{\theta_l}{2}}|) +r_0f'(2r_0 \sin{ |\frac{\theta_l}{2}|})|\sin{\frac{\theta_l}{2}}| \right] \left(1-e^{i(m+1)\theta_l} \right),
\end{equation}

\begin{equation}
    I_2(m):=\sum_{l=1}^{N-1} r_0f'(2r_0|\sin{\frac{\theta_l}{2}}|)|\sin{(\frac{\theta_l}{2})}| \left(e^{i m\theta_l} -e^{i\theta_l}\right).
\end{equation}

\end{proposition}

\begin{proposition}\label{prop2}
There exists an $N$-spot traveling ring solution  to the system \cref{eqhh} with $p_k=v_0 t + r_0e^{i\frac{2k\pi}{N}}$. The radius $r_0$ and velocity $v_0$ are given by 
\begin{subequations}
\begin{align}
        0&=F(r_0),\\
    M_1&=M_2|v_0|^2,
\end{align}
\end{subequations}
whose linear stability is determined by the the eigenvalues of of the following $4\times 4$ matrix for all $m=1,\cdots,N$
 \begin{equation}
     M(m)= 
   \begin{pmatrix}
 -I_1(m)  &  -I_2(m) & 1 & 0 \\
   -I_2(m)& -I_1(-m)  & 0&1\\
  -k_3 I_1(m)  &  -k_3 I_2(m)&-M_1 &-M_1\\
   -k_3 I_2(m)& -k_3 I_1(-m) &-M_1 &-M_1    \end{pmatrix}.
 \end{equation}

\end{proposition}

\begin{proposition}\label{prop3}
There exist an $N$-spot rotating ring solution with $p_k=r_0e^{i(\omega_0 t+\frac{2k\pi}{N})}$ to the system \cref{eqhh} if the following system has a root 
\begin{subequations}\label{omegaandr}
\begin{align}
        \omega_0^2 &=k_3F(r_0) - F^2(r_0),\\
    M_1&=(1+M_2 k_3 r^2_0) F(r_0).
\end{align}
\end{subequations}
Considering the rotating ring solution to the system \cref{eqhh} with radius $r_0$ and frequency $\omega_0$ given by \cref{omegaandr},  we define
\begin{equation}
    H_1=M_1 -2M_2 \left[ \omega_0^2 r_0^2 + r_0^2 F^2(r_0) \right] ,
\end{equation}

\begin{equation}
    H_2=M_2\left[ i\omega_0 r_0 +r_0 F(r_0) \right]^2,
\end{equation}

 \begin{equation}
     M(m)= 
   \begin{pmatrix}
 -I_1(m) -i\omega_0 &  -I_2(m) & 1 & 0 \\
   -I_2(m)& -I_1(-m) +i\omega_0 & 0&1\\
  -k_3 I_1(m)  &  -k_3 I_2(m)&H_1-i\omega_0&-H_2\\
   -k_3 I_2(m)& -k_3 I_1(-m) &-\bar{H}_2&H_1+i\omega_0
    \end{pmatrix},
 \end{equation}
then the rotating $N$-ring is linear stable if the eigenvalues of $M(m)$ have non-positive real parts for all $m=1\cdots N$, 
\end{proposition}

\begin{remark}
   Since $I_1(N-m)=I_1(-m)$ and $I_2(N-m)=I_2(-m)=I_2(m)$,  we only need to compute the eigenvalues of these matrices for $m=0,1,\ldots, \lfloor N \rfloor /2+1 $ to determine the stability.  
\end{remark}

The paper is organized as follows. In  \cref{sec:2}, a reduction from the PDE system to a finite-dimensional ODE systems as described in \cref{lem:1} and \cref{lem:2} is informally presented to extract the nature of the dynamics below and near the drift bifurcation point. This gives rise to a dynamical system with either one phase space dimension or two phase space dimensions per spot and spatial dimension.  In  \cref{sec:3}, we carefully study the reduced ODE system at $\tau<\tau_c$.  The existence and stability of stationary $N$-spot ring solutions are established.  In  \cref{sec:4}, we investigate the reduced ODE system at $\tau \sim \tau_c$. Following the strategy of the stability analysis developed in \cite{kolokolnikov2011stability}, we construct traveling and rotating $N$-spot ring solutions.  The stability of these moving $N$-spot rings is simplified to study the stability of multiple $4\times 4$ matrices. In  \cref{sec:5}, numerical simulations of the PDE and ODE are compared to validate our results.   The conclusion and outlook are presented in  \cref{sec:6}.

\section{Reduced models for the dynamics of N spots}\label{sec:2}
In this section, we briefly derive the reduced ODEs in  \cref{lem:1} and  \cref{lem:2} by center manifold reduction combined with multi-time scale analysis. For a detailed derivation, we refer the reader to \cite{bode2002interaction}. Rigorous results can be found in  \cite{ei2006interacting} for a general system.

\begin{figure}
    \centering
    \includegraphics[width=\textwidth]{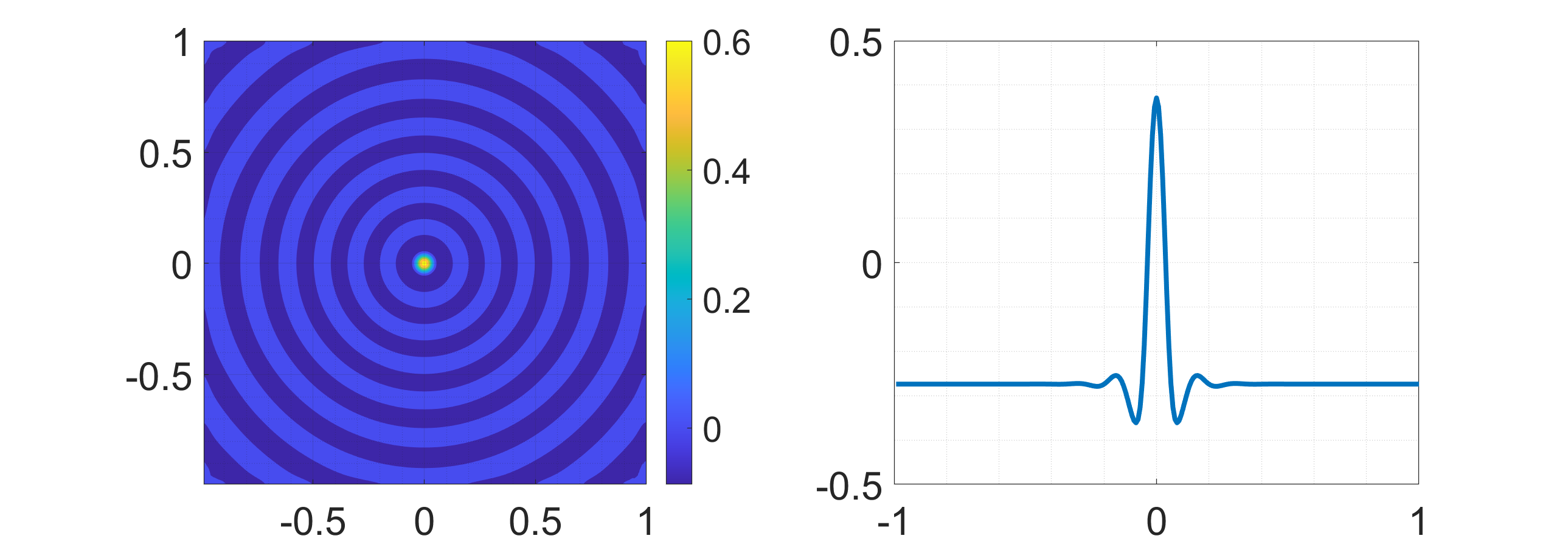}
    \caption{The profile of a single spot steady state for the PDE. Left: the activator's $2$-D contour after subtracting the homogeneous state; Right: the intersection of the activator $u(x)$  at the line $y=0$.}
    \label{fig:spotprofile}
\end{figure}

 For succinctness, we will use $p_k=[x_k,y_k],~q_k=[\xi_k,\eta_k],~\br=[x,y]$ to identify a two-dimensional vector. We represent the inner product of $u$ and $v$ as $ \langle u , v \rangle := \iint_{\mathbb{R}^2} uv \,dx\,dy$.  We assume that the system \cref{rd} admits a stable spot solution in the polar coordinate, denoted as
\begin{equation}
    \mathbf{U}_s(\rho)=[u_c,~u_c]^\intercal+ [ u_s(\rho) ,~u_s(\rho)] ^\intercal \quad \text{with}~~\lim_{\rho \rightarrow \infty} \mathbf{U}_s(\rho) = [u_c,~u_c]^\intercal
\end{equation}
where $u_c$ is a constant corresponding to the homogeneous solution of system \cref{rd} and $u_s$ decays exponentially to $0$ with oscillatory tails, see \cref{fig:spotprofile} for the profile of a single spot.
For the system \cref{rd} with periodic boundary conditions, any translation of the spot is still a solution. For convenience, we define
\begin{equation}
    S_k:=[S_k^u,S_k^v]^\intercal=[u_s(|\br-p_k|),~u_s(|\br-p_k|) ]^\intercal.
\end{equation}
The linearization of the system \cref{rd} for the spot located at $p_k$  gives rise to the operator $\mathcal{L}$
\begin{equation}
    \mathcal{L} :=
    \begin{pmatrix}
    D_u \Delta + k_1 - 3 \left({u}_c+ S_k \right)^2-k_4 \mathcal{G}^{-1} & -k_3\\
    \frac{1}{\tau}& -\frac{1}{\tau}  
    \end{pmatrix}.
\end{equation}
Due to the translation invariance, the operator $\mathcal{L}$ has a eigenvalue $0$ whose corresponding eigenvectors are translational modes:
\begin{equation}\label{Sk}
S_{k,x}=\frac{\partial S_k }{\partial x}, \quad S_{k,y}=\frac{\partial S_k }{\partial y},
\end{equation}
which satisfy 
\begin{equation}\label{def1}
\mathcal{L} S_{k,x}  = 0,~\mathcal{L} S_{k,y}  = 0.
\end{equation}
Similar properties also hold for the adjoint operator $\mathcal{L}^\dagger$ 
\begin{equation}
    \mathcal{L}^\dagger: =
    \begin{pmatrix}
    D_u \Delta + k_1 - 3 \left({u}_c+ S_k \right)^2-k_4 \mathcal{G}^{-1} & \frac{1}{\tau} \\
   -k_3 & -\frac{1}{\tau}  
    \end{pmatrix}.
\end{equation}
There exist eigenvectors
\begin{equation}\label{Skstar}
S_{k,x}^*=\frac{\partial S^*_k}{\partial x},\quad S_{k,y}^* =\frac{\partial S^*_k}{\partial y} \quad \text{with}\quad S^*_{k}:=[{S^*_{k}}^u,{S^*_{k}}^v]^\intercal=[ u_s(|\br-p_k|),~-\tau k_3 u_s(|\br-p_k|)]^\intercal
\end{equation}
such that 
\begin{equation}
 \mathcal{L}^{\dagger} S^*_{k,x}  = 0,\quad\mathcal{L}^\dagger S^*_{k,y}  = 0.
\end{equation}
When $\tau<\tau_c:=\frac{1}{k_3}$, $\mathcal{L}$ has no eigenvalue with positive real part because a single spot is stable by assumption. While at $\tau=\tau_c$, the eigenvalue $0$ is degenerated, there exist generalized eigenvectors to $\mathcal{L}$
 \begin{equation}\label{Psik}
      \Psi_{k,x}= \frac{\partial \Psi_k }{\partial x} ,\quad \Psi_{k,y}=\frac{\partial \Psi_k}{\partial y} \quad \text{with} \quad \Psi_k:=[\Psi_k^u,\Psi_k^v]^\intercal=[0,~\frac{1}{k_3}  u_s(|\br-p_k|) ]^\intercal
 \end{equation}
such that
\begin{gather}\label{def2}
\mathcal{L} \Psi_{k,x}  = -S_{k,x},\quad\mathcal{L} \Psi_{k,y}  = -S_{k,y}.
\end{gather}
For $\mathcal{L}^\dagger$, the generalized eigenvectors are
\begin{equation}\label{Psikstar}
    \Psi_{k,x}^*=\frac{\partial \Psi^*_k}{\partial x}, \quad \Psi_{k,y}^* =\frac{\partial \Psi^*_k}{\partial y}\quad\text{with}\quad \Psi^*_{k}:=[{\Psi^*_{k}}^u,{\Psi^*_{k}}^v]^\intercal =[\frac{1}{k_3} u_s(|\br-p_k|),~0]^\intercal,
\end{equation}
such that
\begin{gather}
    \mathcal{L}^{\dagger} \Psi^*_{k,x}  = -S^*_{k,x},\quad\mathcal{L}^\dagger \Psi^*_{k,y}  = -S^*_{k,y}.
\end{gather}
Thus, to describe the solution of the linear system 
\begin{equation}
    \mathcal{L} \tilde{\mathbf{U}} = P,
\end{equation}
associated with the perturbation of  a spot solution, 
we need to add the generalized eigenvectors to the eigenvector expansion closed to $\tau=\tau_c$, resulting in the expansion in \cref{expansion}.

As the spots are localized and decay exponentially, their superposition is a good approximation to the exact solution when all the distances between these spots are large enough, see \cite{elphick1990patterns}. We consider a superposition of $N$ spots at different positions $p_k,~k=1,\ldots,N$, denoted as
\begin{equation}
    \mathbf{U}_s(r)=[{u}_c,{u}_c]^\intercal + \sum_{k=1}^N S_k.
\end{equation}
   The error of this approximation to the ture solution scales with the shortest distance between spots, as shown by $$\sigma:=\exp{\bigg\{-\frac{cd_{\text{min}}}{\sqrt{D_u}} \bigg\} },$$ where $c$ is a constant related to the spot profile and $d_{\text{min}}$ represents the minimal distance between two neighboring spots. Near the center of the $k$-th spot, the influence from other spots can be interpreted as a perturbation. With this in mind, we proceed to investigate the slow dynamics of multiple spots using perturbation techniques. As the system undergoes a bifurcation at $\tau_c$,  the discussion is split into the following two cases according to the value of $\tau$:

 \textbullet\  When $\tau<\tau_c$, we expand the solution for the dynamics of \cref{rd} as 
 \begin{equation}\label{expansion0}
    \mathbf{U}= \mathbf{U}_s +  \sigma \mathbf{U}_1  +\sigma^2  \mathbf{U}_2+\cdots .
\end{equation}
We assume that each spot moves at a slow time scale $T=\sigma t$, i.e. $p_k=p_k(T)$. Substituting \cref{expansion0} into \cref{rd}, we obtain the corresponding system in power of $\sigma$ near the center of $k$-th spot:
 \begin{fleqn}
 \begin{equation*}
     \mathcal{O}(\sigma): 
 \end{equation*}
 \end{fleqn}
\begin{equation}\label{bbifur}
       - \frac{\partial p_k }{ \partial T} \cdot \nabla S_k=  \mathcal{L} \mathbf{U}_1 +  \begin{pmatrix}
       \frac{1}{\sigma}\left(3 {S^u_k}^2+6 u_c {S_k^u}  \right)  \sum_{l\neq k} {S_l^u} \\0
        \end{pmatrix}.
\end{equation}
Using \cref{Sk} and \cref{Skstar}, taking the inner product of \cref{bbifur} with  $S_{k,x}^*$ and $S_{k,y}^*$, and returning to the original time $t$ yields
 \begin{equation}\label{ordinary}
      \frac{\partial p_k}{\partial t}=-\frac{1}{1-\tau k_3} \sum_{j\neq k} (p_k-p_j) f(|p_k-p_j|),
 \end{equation}
   where
\begin{equation}\label{def_f}
 f(d): = \frac{\iint_\Omega u_{s,x} \left( 3 {u_s}^2+6 u_c {u_s}  \right) u_s( |x-d|,0 ) \,dx\,dy  }{ d \iint_\Omega u_{s,x}^2 \,dx\,dy } .
\end{equation}

For the parameters given in the caption of \cref{N-spot}, the interaction function can be approximated by the following fitting function.
\begin{equation}\label{fitting0}
f(d)= \frac{6.87\times 10^{-4}}{d^{\frac{3}{2}}} e^{-15.7d} \cos{ 43.15(d-0.199) },\quad \text{for}~d>d_b\sim 0.12,
\end{equation}
where $d_b$ represents the core radius of the spot, below which two spots may coalesce into one spot or get annihilated. We note that this fitting formula's analytical form is selected in accordance with the far field behaviour of $u_s$.  The left figure in \cref{fig:M1} gives the comparison between numerically computed $f$ and its fitting approximation. Thus, the movement of the $k$-th spot can be predicted by \cref{ordinary} when other spots are far away with a distance greater that some $d_b$. 
   
 \textbullet\ In the neighborhood of the drift bifurcation point  $\tau=\tau_c+\hat{\tau}\varepsilon^2$ with small parameter $\varepsilon=\sqrt{\sigma} \ll 1$, an appropriate approximate solution for the dynamics of \cref{rd}  is

\begin{equation}\label{expansion}
    \mathbf{U}= \mathbf{U}_s + \varepsilon  \sum_{k=1}^{N} \left( \tilde{q}_k \cdot  \nabla \Psi_k \right)+ \varepsilon^2 \mathbf{U}_2  +\varepsilon^3  \mathbf{U}_3+\cdots .
\end{equation}
We also introduce new time-scales $T_j=\varepsilon^j t$ and assume
\begin{equation}
    p_k=p_k(T_1,T_2,T_3),\quad \tilde{q}_k=\tilde{q}_k(T_1,T_2,T_3).
\end{equation}
To get a unique decomposition, we demand
\begin{equation}
    \langle S_{k,x} , \mathbf{U}_j \rangle=0 \quad \langle S_{k,y} , U_j \rangle=0 \quad
     \langle \Psi_{k,x} , \mathbf{U}_j \rangle =0\quad
      \langle \Psi_{k,y} , \mathbf{U}_j\rangle =0 \quad \text{for}~k=1\ldots N,~\text{and}~j>1,
\end{equation}
where we have defined $\langle \mathbf{f} , \mathbf{g} \rangle=   \iint_{\Omega} f^u g^u+f^v g^v~ \text{d}\mathbf{r}$ for $\mathbf{f}:=[f^u,f^v]^\intercal$ and $\mathbf{g}:=[g^u,g^v]^\intercal$.

The first and second terms of the expansion \cref{expansion} encode the information about spot locations and velocities. They will be balanced in the lowest order of the series. It is worth noting that spot-spot interaction close to the center of the $k$-th spot will appear in the order of $\varepsilon^2$. Substituting \cref{expansion} into \cref{rd}, we obtain the corresponding system in each power of $\varepsilon$ near the center of the $k$-th spot.

\begin{fleqn}
\begin{equation*}
    \mathcal{O}(\varepsilon):
\end{equation*}
\end{fleqn}
\begin{equation}
- \frac{\partial p_k }{ \partial T_1} \cdot \nabla S_k= \tilde{q}_k \cdot \mathcal{L} \nabla \Psi_{k}.
\end{equation}

\begin{fleqn}
\begin{equation*}
    \mathcal{O}(\varepsilon^2):
\end{equation*}    
\end{fleqn}
\begin{multline}
-\frac{\partial p_k }{ \partial T_2} \cdot \nabla S_k    -\nabla (\tilde{q}_k \cdot \nabla \Psi_k) \cdot \frac{\partial p_k }{ \partial T_1} + \frac{\partial \tilde{q}_k}{ \partial T_1}  \cdot \nabla \Psi_k   =  \\
\mathcal{L} \mathbf{U}_2 +      
        \begin{pmatrix}
       \frac{1}{\varepsilon^2}\left(3 {S^u_k}^2+6 u_c {S_k^u}  \right)  \sum_{l\neq k} {S_l^u} \\0
        \end{pmatrix} + \begin{pmatrix}
     \left(3 S_k^u+3 u_c  \right) (\tilde{q}_k \cdot \nabla {\Psi_k^u})^2   \\0
        \end{pmatrix} 
\end{multline}

\begin{fleqn}
\begin{equation*}
    \mathcal{O}(\varepsilon^3):
\end{equation*} 
\end{fleqn}
\begin{multline}
         -\frac{\partial p_k }{ \partial T_3} \cdot \nabla S_k  -\nabla (\tilde{q}_k \cdot \nabla \Psi_k) \cdot \frac{\partial p_k }{ \partial T_2} + \frac{\partial \tilde{q}_k}{ \partial T_2}  \cdot \nabla \Psi_k  + \frac{\partial U_2 }{ \partial T_1}  -  \frac{\hat{\tau}}{\tau_c} \frac{\partial p_k }{ \partial T_1} \cdot \nabla S_k -  \frac{1}{\varepsilon^2}\sum_{l\neq k}\frac{\partial p_l}{\partial T_1} \cdot \nabla S_l   = \\ \mathcal{L} \mathbf{U}_3  + \frac{1}{\varepsilon^2} \sum_{l\neq k}  \mathcal{L} ( \tilde{q}_l \cdot \nabla \Psi_l) +
         \begin{pmatrix}
      (\tilde{q}_k \cdot \nabla {\Psi_k^u})^3   \\0
        \end{pmatrix} +  \begin{pmatrix}
     \left(3 S_k^u+3 u_c  \right) (\tilde{q}_k \cdot \nabla {\Psi_k^u}) U_1^u   \\0
        \end{pmatrix}  
\end{multline}
Using \cref{Sk},~\cref{Skstar},~\cref{Psik} and \cref{Psikstar}, taking the inner product with $S^{*}_{k,x},S^{*}_{k,y},\Psi^{*}_{k,x}$ and $\Psi^{*}_{k,y}$ in each powers of $\varepsilon$ yield

\begin{subequations}\label{herachy}
\begin{align}
    \frac{\partial p_k}{\partial T_1}&= \tilde{q_k} ,\\
     \frac{\partial p_k}{\partial T_2}&=-\frac{1}{\varepsilon^2}\sum_{j\neq k} (p_k-p_j) f(|p_k-p_j|),\\
    \frac{\partial \tilde{q}_k}{\partial T_1}&= - \frac{k_3}{\varepsilon^2} \sum_{j\neq k} (p_k-p_j) f(|p_k-p_j|),\\
     \frac{\partial \tilde{q}_k}{\partial T_2}
   & =k_3^2 \hat{\tau} \tilde{q}_k - \frac{Q}{k_3} \tilde{q}_k |\tilde{q}_k|^2,
\end{align}
\end{subequations}
where $f$ is defined in \cref{def_f} and
\begin{equation}\label{numf}
Q= \frac{ \iint u_{s,xx}^2\,dx\,dy   }{\iint u_{s,x}^2 \,dx\,dy }.
\end{equation}
Note that 
\begin{equation}\label{ddt}
 \dot{p}_k:=\frac{d p_k}{d t}= \varepsilon \frac{\partial p_k}{\partial T_1} + \varepsilon^2 \frac{\partial p_k}{\partial T_2}+\cdots;\quad  \dot{\tilde{q}}_k:= \frac{d \tilde{q_k}}{d t}= \varepsilon \frac{\partial \tilde{q_k}}{\partial T_1} + \varepsilon^2 \frac{\partial \tilde{q_k}}{\partial T_2}+\cdots.
\end{equation}
Substituting \cref{herachy} into \cref{ddt}, neglecting higher order terms and returning to the original variable without the small parameter $\varepsilon$,  we obtain the ODE system for the $p_k$ and $q_k:=\varepsilon \tilde{q}_k$:
\begin{subequations}\label{eqhhode}
\begin{equation}\label{eqhh1ode}
    \dot{p}_k=q_k-\sum_{j\neq k} (p_k-p_j) f(|p_k-p_j|) ,
\end{equation}
\begin{equation}\label{eqhh2ode}
    \dot{q}_k=M_1 q_k - M_2 q_k |q_k|^2 -k_3\sum_{j\neq k} (p_k-p_j) f(|p_k-p_j|) ,
\end{equation}
\end{subequations}
where 
\begin{gather}
    M_1=k_3^2(\tau-\frac{1}{k_3}),~M_2= \frac{Q}{k_3}.
\end{gather}
In this way, the PDE system \cref{rd} near the bifurcation point is reduced to a $4N$-dimensional ODE system \cref{eqhhode} that can be recognized as the normal form of the drift bifurcation. The reduced description,
the ODE system \cref{eqhhode}, provides a powerful tool for quickly and conveniently exploring the motion of multiple spots.

The analysis in the paper's reminder is based on the reduced ODE systems. We emphasize that the reduced systems consider merely the effect of the translation mode. Other modes that may cause the dramatic change of a spot profile occur at a fast time scale and are assumed to be unexcited. Scenarios involving spot creation and destruction are not covered by the ODE. Also, we note that the reduced systems are valid only when the spot-spot distance is large and the system is near the drift bifurcation, because under these conditions, one spot is considered as a small perturbation to the other spot.

\section{Stationary N-spot rings and their stability} \label{sec:3}
In this section, we show the existence of  $N$-spot ring solutions to the system \cref{ordinary}, the stability of which is determined by $N$ matrices of $2\times 2$ size.  \cref{prop1} is a direct result of the analysis.  Throughout the remainder of the paper, we will identify $e^{i\theta} := (\cos{\theta}, \sin{\theta})$ to the corresponding complex number interchangeably when referring to ring solutions and $r_0$ to a scalar referred to as the radius of the ring.

\subsection{ N-spot ring}
\begin{figure}
    \centering
    \includegraphics[width=0.48\textwidth]{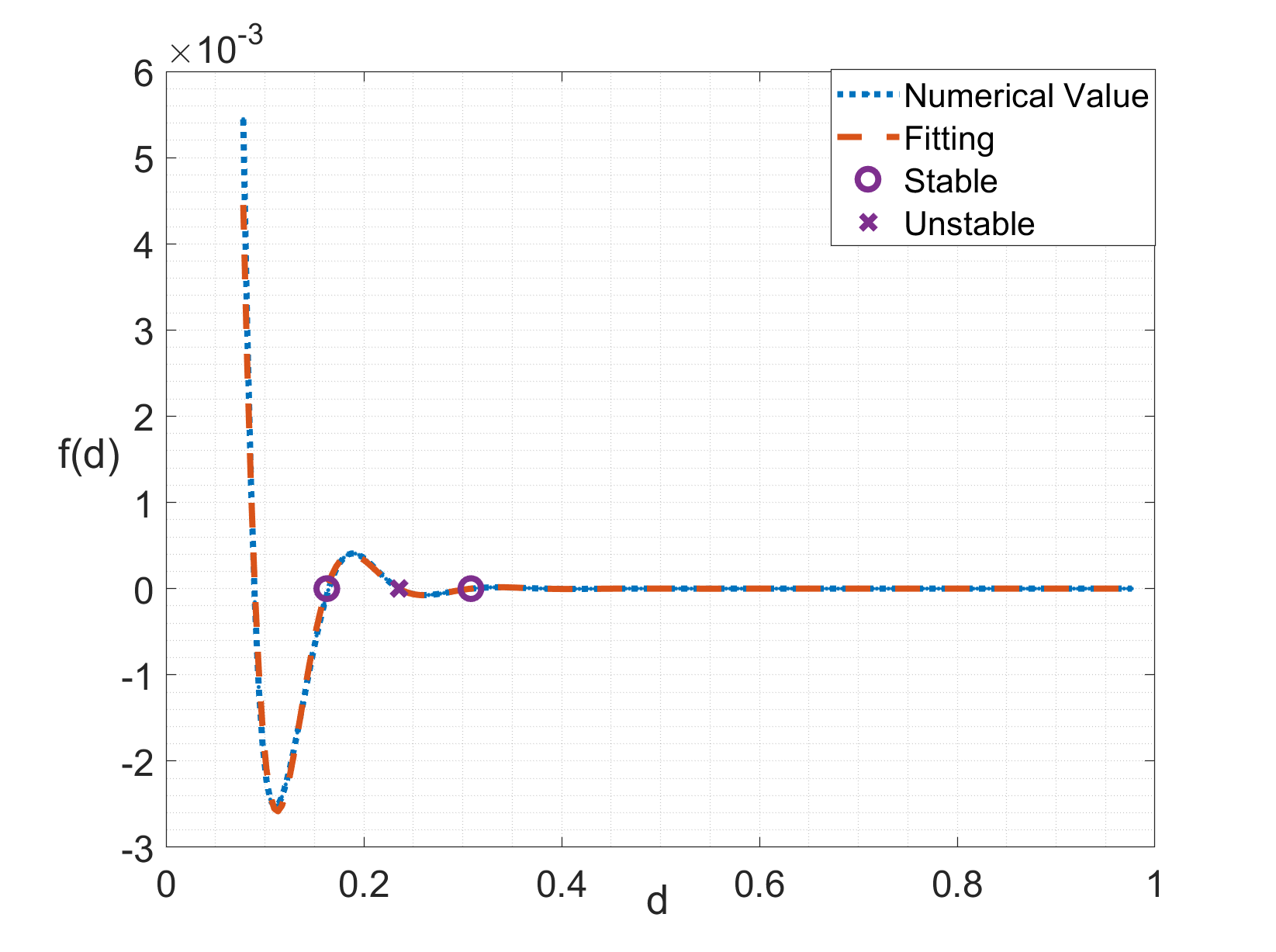}
    \includegraphics[width=0.48\textwidth]{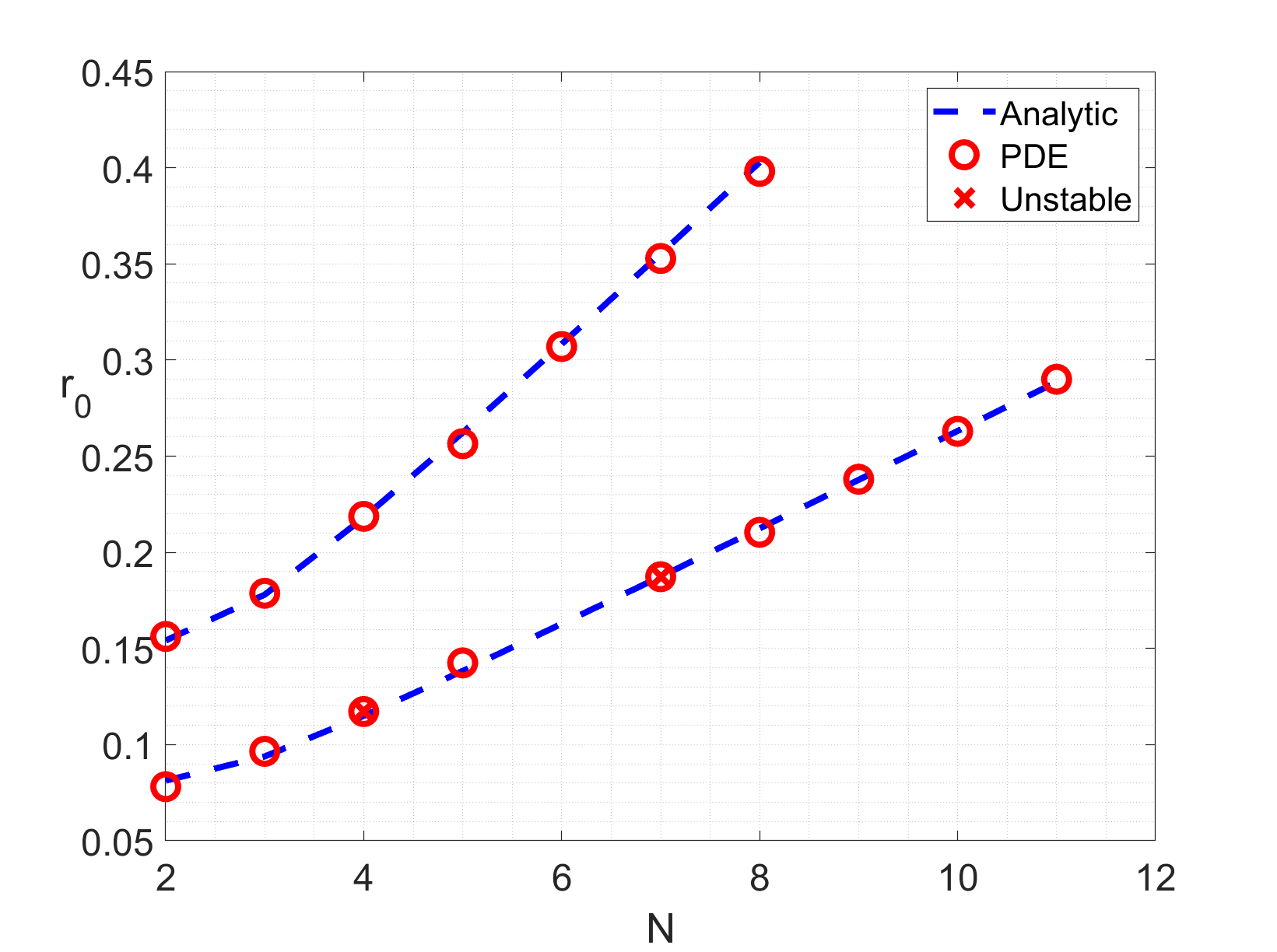}
    \caption{Left: $f(d)$ obtained numerically by \cref{def_f} and the fitting \cref{fitting0}. The numerical value is represented by the dot line, and the fitting value is represented by the red dashed line. The circles indicate the stable distance for a two-spot ring, whereas the crosses mark the unstable distance. Right: The radius of an N-spot ring as a function of the spot count. The radius computed via PDE simulation is denoted by red circles, whereas the analytical value predicted by \cref{r0} is indicated by a blue dashed line. The radius for $N=6$ is omitted from PDE simulation since structural instability is triggered at the center.  A cross is superimposed on the circle to indicate the unstable ring configuration. }
    \label{fig:M1}
\end{figure}
 \begin{figure}
     \centering
     \includegraphics[width=\textwidth]{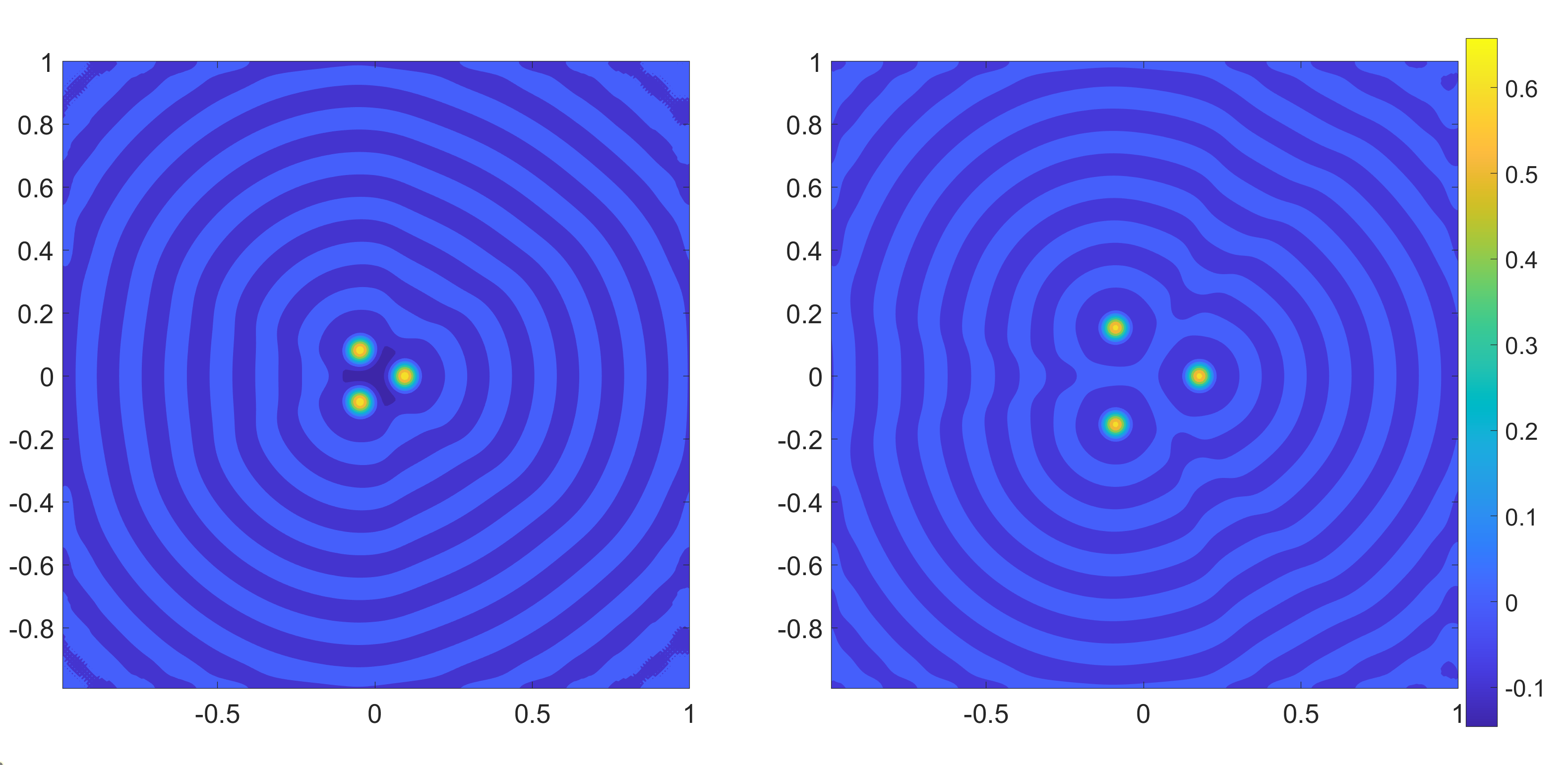}
     \caption{Contours of the activator distribution after subtracting the homogeneous state for three-spot rings with the first and second binding radii. The parameters are given in the caption of \cref{N-spot}. Left: the radius is approximately $0.0965$, slightly above $0.0939$, the first attractive zero of $F(r_0)$; Right: the radius is approximately $0.1788$, slightly above $0.1780$, the second attractive zero of $F(r_0)$.  }
     \label{fig:3spots}
 \end{figure}
 
In this subsection, we construct an $N$-spot ring solution to \cref{ordinary}. In particular, we seek a solution with the form of a ring as follows:
 \begin{equation} \label{assump0}
     p_k=r_0 e^{i\theta_k},\quad\text{where }\theta_k=\frac{2\pi k}{N},\quad k=1,\ldots, N.
 \end{equation}
 The equilibrium point of \cref{ordinary} then satisfies
\begin{equation}\label{ss}
     \sum_{j\neq k} (p_k-p_j) f(|p_k-p_j|)=0.
\end{equation}
 Using the identity
 \begin{equation}\label{identity}
      (p_k-p_{k+l}) f(|p_k-p_{k+l}|)+ (p_k-p_{k+N-l}) f(|p_k-p_{k+N-l}|) 
=   4r_0\sin^2{\frac{\theta_l}{2}} f( 2r_0|\sin{\frac{\theta_l}{2}}|) e^{i \theta_k},
 \end{equation}
and substituting \cref{assump0} back into \cref{ss} yields
 \begin{equation}\label{ssra}
      r_0  F(r_0) e^{i\theta_k}=0,
 \end{equation} 
 where
 \begin{equation}\label{Fr0}
    F(r_0):=\sum_{l=1}^{N-1} (1-e^{i\theta_l })  f( 2r_0|\sin{\frac{\theta_l}{2}}|)=\sum_{j=1}^n s_j  f(c_j r_0),
\end{equation}
with
\begin{equation}
\begin{cases}
    s_j=4\sin^2{\frac{\theta_j}{2} },\quad &c_j=2|\sin{\frac{\theta_j}{2} }|,\quad \text{for }j\neq n,~ \text{when }N=2n+1, \\
    s_j=2,\quad &c_j=2,\quad  \text{for }j= n,~ \text{when }N=2n.
\end{cases}
\end{equation}
As $f(r_0)$ is small and decays exponentially, $f(c_j r_0) \ll f(c_1 r_0) $ for $j\geq2$. Thus $F(r_0)$ can be further approximated by the first component in the series \cref{Fr0},
\begin{equation}\label{Fappro}
    F(r_0) \sim 4 \sin^2{\frac{\pi}{N}} f(2r_0| \sin{\frac{\pi}{N}}| ) \quad \text{when}~ N>2.
\end{equation}
Then the solution to \cref{ssra} can be solved as
\begin{equation}\label{r0}
    r_0\sim \frac{d_c}{ 2| \sin{\frac{\pi}{N}}|},\quad \text{with} ~ d_c ~\text{satisfies}~ f(d_c)=0.
\end{equation}
Due to the oscillatory behavior of $f(d)$ around zero,  $d_c$ can be some discrete values. For the given parameters in {\cref{N-spot}}, these values are $d_{c}=0.1626,0.2354,0.3083, \cdots$, which are shown as circles and crosses in the right figure of  \cref{fig:M1}.  For large $N$, $\sin{\frac{\pi}{N}} \sim \frac{\pi}{N}$, thus $r_0 \sim \frac{Nd_c}{2\pi}$. We will use these approximations as our initial guesses to construct an $N$-spot ring in the numerical simulations.

\subsection{Stability of N-spot rings}
In this subsection, we  analyze the linear stability of the ring solution with radius $r_0$ given by \cref{ssra}. We begin by introducing the perturbations to a ring of $N$ spots in the following form
\begin{equation}\label{perbring0}
     p_k= r_0 e^{i\theta_k}(1 + \phi_k)
\end{equation}
with $\phi_k$  such that $|\phi_k|\ll 1$. Let $l=j-k$, we compute
\begin{equation}
    p_k-p_j=  r_0 e^{i\theta_k} \left(1+\phi_k - e^{i\theta_l}(1+\phi_j)  \right),
\end{equation}
then
\begin{equation}
\begin{aligned}
   |p_k-p_j| &=  r_0 |1-e^{i\theta_l }| +\frac{ r_0}{2 |1-e^{i\theta_l}|}    \left[(1-e^{-i\theta_l})(\phi_k+\bar{ \phi}_j)+ (1-e^{i\theta_l})(\bar{\phi}_k+\phi_j)\right] + h.o.t,\\
   &= 2r_0|\sin{\frac{\theta_l}{2}}|+\frac{ r_0}{4|\sin{\frac{\theta_l}{2}}|} \left[(1-e^{-i\theta_l})(\phi_k+\bar{ \phi}_j)+ (1-e^{i\theta_l})(\bar{\phi}_k+\phi_j)\right] + h.o.t.
\end{aligned}
\end{equation}
Substituting \cref{perbring0} into \cref{ordinary} and neglecting higher order terms leads to the following system:
\begin{subequations}\label{eigenN0}
\begin{multline}
    \dot{\phi}_k = - \frac{1}{1-\tau k_3} \sum_{j\neq k}   (\phi_k- e^{i\theta_l} \phi_j)  f(2r_0|\sin{\frac{\theta_l}{2}}|)\\
   -\frac{1}{1-\tau k_3}\sum_{j\neq k}   \frac{ r_0}{4|\sin{\frac{\theta_l}{2}}|} \left[(1-e^{-i\theta_l})(\phi_k+\bar{ \phi}_j)+ (1-e^{i\theta_l})(\bar{\phi}_k+\phi_j)\right] f'(2r_0|\sin{\frac{\theta_l}{2}}|)(1-e^{i\theta_l})
\end{multline}
\end{subequations}
Using the identity
\begin{equation}
    (1-e^{i\theta_l})^2=-4\sin^2{\frac{\theta_l}{2}}e^{i\theta_l};\quad (1-e^{i\theta_l})(1-e^{-i\theta_l})=4\sin^2{\frac{\theta_l}{2}},
\end{equation}
we obtain
\begin{subequations}\label{eigenNN0}
\begin{equation}
    \dot{\phi}_k = -\frac{1}{1-\tau k_3} \sum_{j\neq k} \left( G_1(\frac{\theta_l}{2}) (\phi_k-e^{i\theta_l} \phi_j)+ G_2(\frac{\theta_l}{2}) ( \bar{\phi}_j-e^{i\theta_l} \bar{\phi}_k  ) \right)
\end{equation}
\end{subequations}
with 
\begin{equation}
\begin{aligned}
    G_1(\theta)&= f(2r_0|\sin{\theta}|) +r_0f'(2r_0|\sin{\theta|})|\sin{\theta}|, \\
    G_2(\theta)&=r_0f'(2r_0|\sin{\theta|})|\sin{\theta}|.
\end{aligned}
\end{equation}
Assuming that $\phi_k$ satisfies the following relation
\begin{equation}\label{ansatz0}
    \phi_k=\xi_+ e^{i m \theta_k}+\xi_- e^{-i m \theta_k} \quad m = 1,\ldots, N,
\end{equation}
then we can write $\phi_j$ as
\begin{equation}\label{rewrite0}
    \phi_j = \xi_+ e^{i m \theta_k} e^{i m \theta_l } +\xi_- e^{-i m \theta_k} e^{-i m \theta_l}, \quad m = 1,\ldots, N.
\end{equation}
Substituting \cref{ansatz0} and \cref{rewrite0} into \cref{eigenN0} and collecting like terms in $e^{im\theta_k},e^{-im\theta_k}$ leads to the system
\begin{subequations}\label{Neigeneq0}
\begin{equation}
  \dot{\xi}_+ = -\frac{1}{1-\tau k_3} \xi_+ \sum_{j\neq k}  G_1(\frac{\theta_l}{2}) \left(1-e^{i(m+1)\theta_l} \right)- \frac{1}{1-\tau k_3}\bar{\xi}_- \sum_{j\neq k} G_2(\frac{\theta_l}{2}) \left(e^{im\theta_l}-e^{i\theta_l} \right) ,
\end{equation}

\begin{equation} \label{xi0-}
    \dot{\xi}_- = - \frac{1}{1-\tau k_3}\xi_- \sum_{j\neq k}    G_1(\frac{\theta_l}{2}) \left(1-e^{i(-m+1)\theta_l} \right) -\frac{1}{1-\tau k_3}\bar{\xi}_+ \sum_{j\neq k}    G_2(\frac{\theta_l}{2}) \left(e^{-im\theta_l}-e^{i\theta_l} \right) 
\end{equation}
\end{subequations}
Note that the sums are independent of $k$ and $j$. We define
\begin{equation}
    I_1(m)=\sum_{l=1}^{N-1} G_1(\frac{\theta_l}{2}) \left(1-e^{i(m+1)\theta_l} \right)=4 \sum_{l=1}^{N/2}G_1( \frac{\pi l}{N} ) \sin^2{ \frac{ (m+1) \pi l}{N} }
\end{equation}

\begin{equation}
    I_2(m)=\sum_{l=1}^{N-1} G_2(\frac{\theta_l}{2}) \left(e^{im\theta_l} -e^{i\theta_l}\right)=4 \sum_{l=1}^{N/2}G_2( \frac{\pi l}{N} ) \left(  \sin^2{ \frac{ \pi l}{N} }- \sin^2{ \frac{ m \pi l}{N} } \right)
\end{equation}
Using these notations and taking a conjugate of \cref{xi0-} yields
\begin{equation} \label{matrix0}
    \begin{pmatrix} \dot{\xi}_+ \\
   \dot{\bar{\xi}}_-
   \end{pmatrix}=\frac{1}{1-\tau k_3}
   \begin{pmatrix}
 -I_1(m) &  -I_2(m)  \\
   -I_2(m)& -I_1(-m) 
    \end{pmatrix}
      \begin{pmatrix} \xi_+ \\
   \bar{\xi}_-
   \end{pmatrix}
\end{equation}
Let  $\begin{pmatrix} \xi_+ \\
   \bar{\xi}_-
   \end{pmatrix}=e^{\lambda t}
   \begin{pmatrix}
   a_1\\
   a_2
   \end{pmatrix}$, then $\lambda$ is an eigenvalue of the $2\times 2$ matrix in \cref{matrix0}. The eigenvalue $\lambda$ must satisfy:
   \begin{equation}
       \lambda a =\frac{1}{1-\tau k_3} G a,
   \end{equation}
 where 
 \begin{equation}\label{G}
     G:= 
   \begin{pmatrix}
 -I_1(m)  &  -I_2(m)\\
   -I_2(m)& -I_1(-m) 
    \end{pmatrix}.
 \end{equation}
In this way, the study of stability of an $N$-spot ring solution decouples into the study of individual Fourier modes. We conclude the results in  \cref{prop1}.

As the impacts from other non-neighbor spots diminish exponentially according to the distance, they are relatively insignificant in comparison to the neighboring spots. Thus we further approximate the summation of the interaction terms by the interaction terms between a spot and its nearest two neighbor spots when $N$ is large, yielding the following approximations:
\begin{equation}\label{Iappr}
    I_1(m) \sim 4 G_1(\frac{\pi}{N}) \sin^2{ \frac{ (m+1) \pi }{N} },\quad I_2(m) \sim 4 G_2( \frac{\pi}{N} ) \left(  \sin^2{ \frac{ \pi }{N} }- \sin^2{ \frac{ m \pi }{N} } \right)
\end{equation}
 Note that  
 \begin{equation}\label{Gappr}
\begin{aligned}
    G_1( \frac{\pi}{N} )&= f(2r_0|\sin{\frac{\pi}{N}}|) +r_0f'(2r_0 |\sin{ \frac{\pi}{N} }|)|\sin{\frac{\pi}{N}|} \sim \frac{d_c}{2} f'(d_c) \\
    G_2(\frac{\pi}{N} )&=r_0f'(2r_0|\sin{\frac{\pi}{N}}|)|\sin{\frac{\pi}{N}}|\sim \frac{d_c}{2} f'(d_c) .
\end{aligned}
 \end{equation}
Using \cref{Iappr} and \cref{Gappr}, we obtain
 \begin{equation}\label{G_approx}
     G\sim -2d_c f'(d_c)
   \begin{pmatrix}
 \sin^2{ \frac{ (m+1) \pi }{N} }  &    \sin^2{ \frac{ \pi }{N} }- \sin^2{ \frac{ m \pi }{N} } \\
     \sin^2{ \frac{ \pi }{N} }- \sin^2{ \frac{ m \pi }{N} } & \sin^2{ \frac{ (-m+1) \pi }{N} }
    \end{pmatrix}.
 \end{equation}
The eigenvalues of $G$ can be computed directly as 
\begin{equation}\label{eigenexact}
    \lambda=0,~\text{or}~-2d_c f'(d_c)\left(  \sin^2{ \frac{ (m+1) \pi }{N} } + \sin^2{ \frac{ (-m+1) \pi }{N} }  \right) . 
\end{equation}
Thus we arrived at the following corollary:
\begin{corollary}
Under the assumption that only the nearest spot-spot interaction is considered, the $N$-spot ring with the spot-spot distance $d_c$ is stable when $f'(d_c)>0$.
\end{corollary}

 As $f(d)$ oscillates around zero, the stable and unstable distance appears alternatively. Roots of $f(d)=0$ then can be classified as attractive zeros ($f'(d_c)>0$) and repulsive zeros ($f'(d_c)<0$).   By convention, we denote the $N$-spot ring with the spot-spot distance closed to the $k$-th attractive zero point of $f$ as the $N$-spot ring with the $k$-th binding radius.  Henceforth, we always take $d_c$ as the stable zeros of $f$.  Additionally, we exclude the neutral modes coming from $N$ zero eigenvalues. Therefore our corollary states that all the eigenvalues have non-positive real parts, namely "linearly stable".  The effects coming from higher-order terms are important to judge how those neutral modes behave, however we do not go into the details here, but our conjecture is the following: 
 \begin{conjecture}
 All the N-spot ring patterns are ``nonlinear stable" unless spot-spot distance is the smallest binding distance.
 \end{conjecture}
 The conjecture is confirmed numerically for $N=2,\ldots,8$. In other words, we distinguish linear stable and nonlinear stable patterns. We can make a rigorous statement in the former case, but only a conjecture for the latter case.   In \cref{fig:3spots}, two stable three-spot rings corresponding to the first two binding radii are depicted. We remark that using the full summation of the interaction terms, $N$-spot rings with the first binding radius are unstable when $N=4,7$ under the parameter setting in \cref{N-spot}, which is confirmed with PDE simulation.

\section{Traveling N-spot rings and rotating N-spot rings near the drift bifurcation and their linear stability } \label{sec:4}
\begin{figure}[!htb]
    \centering
    \includegraphics[width=0.9\textwidth]{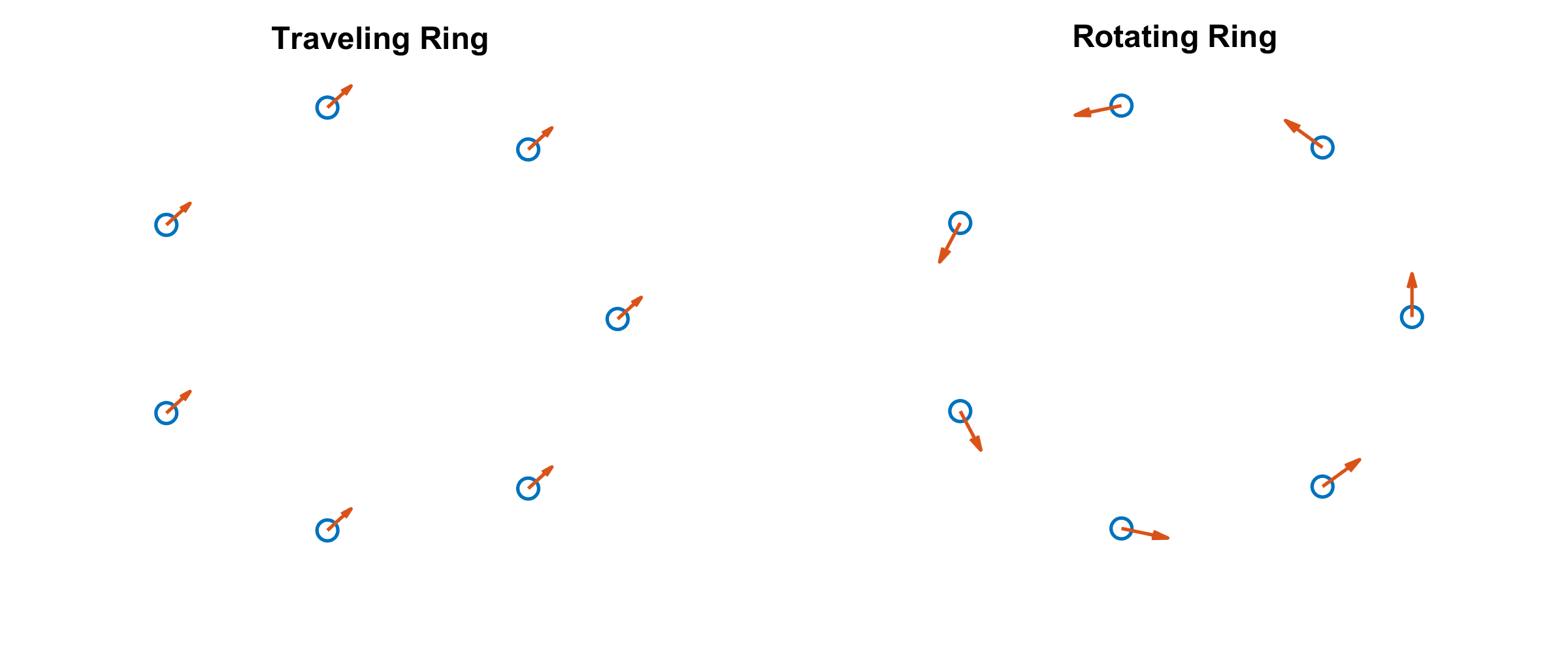}
    \caption{Traveling and rotating rings.}
    \label{fig:illus}
\end{figure}

In this section, we construct two particular solutions  to the reduced ODE system \cref{eqhhode}: traveling $N$-spot rings and rotating $N$-spot rings (see \cref{fig:illus}).  The stability of these two solutions are determined by the eigenvalues of $N$ matrices of $4\times 4$  size.   \cref{prop2} and \cref{prop3} are obtained as the results of the analysis.

\subsection{Traveling N-spot rings}
We start by constructing a traveling ring solution to the reduced ODE system \cref{eqhhode}. By abuse of notation, we seek a particular solution to \cref{eqhhode} with the form of a traveling ring as follows
 \begin{equation} \label{assumpt}
     p_k=v_0 t + r_0 e^{i \theta_k},\quad\text{where }\theta_k=\frac{2\pi k}{N},\quad k=1,\ldots, N,
 \end{equation}
 where $v_0=v_{0,x}+i v_{0,y}$ is the velocity of the ring. 
 Substituting \cref{assumpt} back into \cref{eqhh1}, we obtain
 \begin{equation}\label{q_kt}
     q_k= v_0 + r_0 F(r_0)e^{i\theta_k},
 \end{equation}
 where $F(r_0)$ is defined in \cref{Fr0}. 
 Substituting \cref{q_kt} into \cref{eqhh2} yields
 \begin{equation}\label{eq6-0-t}
   0 = \left(M_1- M_2 |v_0 + r_0 F(r_0) e^{i\theta_k}|^2 \right) \left(v_0 + r_0 F(r_0) e^{i\theta_k} \right)  - k_3 r_0  F(r_0) e^{i\theta_k}. 
 \end{equation}
Equating the imaginary part and the real part of \cref{eq6-0-t}  gives rise to
\begin{align}\label{omega_r_t}
        0&=M_1-M_2|v_0|^2\\
    0&= F(r_0).
\end{align}
From \cref{omega_r_t}, a traveling ring solution has a fixed radius and moves along any direction with a fixed magnitude.

\subsection{Stability of traveling N-spot rings}
In this subsection, we  analyze the stability of the traveling ring solution of radius $r_0$ with velocity $v_0$ given by \cref{omega_r_t}. We first introduce a perturbations to the traveling $N$-spot ring in the following form
\begin{equation}\label{perbringt}
     p_k=v_0 t + r_0 e^{i\theta_k}(1 + \phi_k), \quad    q_k= v_0 + e^{i\theta_k} \psi_k.
\end{equation}
with $\phi_k$ and $\psi_k$ such that $|\phi_k|,|\psi_k|\ll 1$. Let $l=j-k$, we compute
\begin{equation}
\begin{aligned}
    |q_k|^2 & = (v_0+ e^{i\theta_k} \psi_k)  (\bar{v}_0+ e^{-i\theta_k} \bar{\psi_k})\\
                & = |v_0|^2+\bar{v}_0 e^{i\theta_k} \psi_k + v_0 e^{-i\theta_k} \bar{\psi}_k+|\psi_k|^2.
\end{aligned}
     \end{equation}
Substituting \cref{perbringt} into \cref{eqhhode} and neglecting higher order terms leads to the following system:
\begin{subequations}\label{eigenN_t}
\begin{multline}
    \dot{\phi}_k =  \psi_k - \sum_{j\neq k}   (\phi_k- e^{i\theta_l} \phi_j)  f(2r_0|\sin{\frac{\theta_l}{2}}|)\\
   -\sum_{j\neq k}   \frac{ r_0}{4|\sin{\frac{\theta_l}{2}}|} \left[(1-e^{-i\theta_l})(\phi_k+\bar{ \phi}_j)+ (1-e^{i\theta_l})(\bar{\phi}_k+\phi_j)\right] f'(2r_0|\sin{\frac{\theta_l}{2}}|)(1-e^{i\theta_l})
\end{multline}
\begin{multline}
     \dot{\psi}_k =\left[ \left( M_1-2M_2 |v_0|^2\right] \psi_k  - M_2v_0^2  \bar{\psi}_k \right)
    -k_3 \sum_{j\neq k}   (\phi_k- e^{i\theta_l} \phi_j)  f(2r_0|\sin{\frac{\theta_l}{2}}|)\\
    -k_3\sum_{j\neq k}   \frac{ r_0}{4|\sin{\frac{\theta_l}{2}}|} \left[(1-e^{-i\theta_l})(\phi_k+\bar{ \phi}_j)+ (1-e^{i\theta_l})(\bar{\phi}_k+\phi_j)\right] f'(2r_0|\sin{\frac{\theta_l}{2}}|)(1-e^{i\theta_l}).
\end{multline}
\end{subequations}
With the same notations as in  \cref{sec:3}, we write \cref{eigenN_t} as
\begin{subequations}\label{eigenNNt}
\begin{equation}
    \dot{\phi}_k = \psi_k - \sum_{j\neq k} \left( G_1(\frac{\theta_l}{2}) (\phi_k-e^{i\theta_l} \phi_j)+ G_2(\frac{\theta_l}{2}) ( \bar{\phi}_j-e^{i\theta_l} \bar{\phi}_k  ) \right)
\end{equation}
\begin{equation}
     \dot{\psi}_k =\left[ \left(M_1-2M_2|v_0|^2\right) \psi_k  - M_2v_0^2  \bar{\psi}_k \right]
   -k_3 \sum_{j\neq k} \left( G_1(\frac{\theta_l}{2}) (\phi_k-e^{i\theta_l} \phi_j) + G_2(\frac{\theta_l}{2}) ( \bar{\phi}_j-e^{i\theta_l} \bar{\phi}_k  ) \right)
\end{equation}
\end{subequations}
Assuming that $\phi_k$ satisfies the following relation
\begin{equation}\label{ansatz_t}
    \phi_k=\xi_+ e^{i m \theta_k}+\xi_- e^{-i m \theta_k},\quad \psi_k=\eta_+ e^{i m \theta_k}+\eta_- e^{-i m \theta_k}, \quad m = 1,\ldots, N.
\end{equation}
Then we can write $\phi_j$ as
\begin{equation}\label{rewrite_t}
    \phi_j = \xi_+ e^{i m \theta_k} e^{i m \theta_l } +\xi_- e^{-i m \theta_k} e^{-i m \theta_l}, \quad m = 1,\ldots, N.
\end{equation}
Substituting \cref{ansatz_t} and \cref{rewrite_t} into \cref{eigenN_t} and collecting like terms in $e^{im\theta_k},e^{-im\theta_k}$ leads to the system
\begin{subequations}\label{Neigeneq_t}
\begin{equation}
   \dot{\xi}_+ = \eta_+ -  \xi_+ \sum_{j\neq k}  G_1(\frac{\theta_l}{2}) \left(1-e^{i(m+1)\theta_l} \right)- \bar{\xi}_- \sum_{j\neq k} G_2(\frac{\theta_l}{2}) \left(e^{im\theta_l}-e^{i\theta_l} \right) ,
\end{equation}

\begin{equation} \label{xi-t}
     \dot{\xi}_- = \eta_- - \xi_- \sum_{j\neq k}    G_1(\frac{\theta_l}{2}) \left(1-e^{i(-m+1)\theta_l} \right) -\bar{\xi}_+ \sum_{j\neq k}    G_2(\frac{\theta_l}{2}) \left(e^{-im\theta_l}-e^{i\theta_l} \right) 
\end{equation}

\begin{equation}
      \dot{\eta}_+ = \left(  M_1 -2M_2|v_0|^2 \right) \eta_+  -  M_2v_0^2 \bar{\eta}_-  -k_3 \xi_+ \sum_{j\neq k}  G_1(\frac{\theta_l}{2}) \left(1-e^{i(m+1)\theta_l} \right)-  k_3\bar{\xi}_- \sum_{j\neq k} G_2(\frac{\theta_l}{2}) \left(e^{im\theta_l}-e^{i\theta_l} \right) 
\end{equation}

\begin{equation}\label{eta-t}
     \dot{\eta}_- =   \left(M_1 -2M_2 |v_0|^2 \right) \eta_- -  M_2v_0^2\bar{\eta}_+      -k_3 \xi_- \sum_{j\neq k}    G_1(\frac{\theta_l}{2}) \left(1-e^{i(-m+1)\theta_l} \right) -k_3\bar{\xi}_+  \sum_{j\neq k}    G_2(\frac{\theta_l}{2}) \left(e^{-im\theta_l}-e^{i\theta_l} \right) 
\end{equation}
\end{subequations}
Using the notations in  \cref{sec:3} and taking a conjugate of \cref{xi-t} and \cref{eta-t}  yields
\begin{equation} \label{matrixt}
    \begin{pmatrix} \dot{\xi}_+ \\
   \dot{\bar{\xi}}_-\\
   \dot{\eta}_{+}
   \\
   \dot{\bar{\eta}}_{-}
   \end{pmatrix}=
   \begin{pmatrix}
 -I_1(m) &  -I_2(m) & 1 & 0 \\
   -I_2(m)& -I_1(-m)  & 0&1\\
  -k_3 I_1(m)  &  -k_3 I_2(m)&M_1 -2M_2 |v_0|^2&-M_2 v_0^2\\
   -k_3 I_2(m)& -k_3 I_1(-m) &-M_2 \bar{v}_0^2&M_1 -2M_2 |v_0|^2
    \end{pmatrix}
      \begin{pmatrix} \xi_+ \\
   \bar{\xi}_-\\
   \eta_{+}
   \\
   \bar{\eta}_{-}
   \end{pmatrix}
\end{equation}
Let  $\begin{pmatrix} \xi_+ \\
   \bar{\xi}_-\\
   \eta_{+}
   \\
 \bar{\eta}_{-}
   \end{pmatrix}=e^{\lambda t}
   \begin{pmatrix}
   a_1\\
   a_2\\
   a_3\\
   a_4
   \end{pmatrix}$, then $\lambda$ is an eigenvalue of the $4\times 4$ matrix in \cref{matrixt}. The eigenvalue $\lambda$ must satisfy:
   \begin{equation}
       \lambda a = M a
   \end{equation}
 where 
 \begin{equation}\label{M-t}
     M= 
    \begin{pmatrix}
 -I_1(m) &  -I_2(m) & 1 & 0 \\
   -I_2(m)& -I_1(-m)  & 0&1\\
  -k_3 I_1(m)  &  -k_3 I_2(m)& -M_2 |v_0|^2&-M_2 v_0^2\\
   -k_3 I_2(m)& -k_3 I_1(-m) &-M_2 \bar{v}_0^2&-M_2 |v_0|^2
    \end{pmatrix}.
 \end{equation}
Therefore, the stability of a traveling N-spot ring with velocity $(v_{0,x},~v_{0,y})$ is determined by the eigenvalue of $M$ in \cref{M-t} for $m=1\cdots N$.  By rotational invariance, it suffices to only consider the case $v_0=\sqrt{M_1/M_2}$, namely,
 \begin{equation}\label{M-tt}
     M= 
    \begin{pmatrix}
 -I_1(m) &  -I_2(m) & 1 & 0 \\
   -I_2(m)& -I_1(-m)  & 0&1\\
  -k_3 I_1(m)  &  -k_3 I_2(m)& -M_1 &-M_1 \\
   -k_3 I_2(m)& -k_3 I_1(-m) &-M_1 &-M_1 
    \end{pmatrix}.
 \end{equation}
In this way, the study of stability of an traveling $N$-spot ring solution decouples into the study of individual Fourier modes. We conclude the results in \cref{prop2}.

\subsection{Rotating N-spot rings}
In this subsection, we construct rotating ring solution to the reduced ODE system \cref{eqhhode} for further analysis.
By abuse of notation, we seek a solution to \cref{eqhhode} with the form of a rotating ring as follows
 \begin{equation} \label{assump}
     p_k=r_0 e^{i(\omega_0 t +\theta_k)},\quad\text{where }\theta_k=\frac{2\pi k}{N},\quad k=1,\ldots, N.
 \end{equation}
 Substituting \cref{assump} back into \cref{eqhh1} and using the identity
 \begin{equation}
      (p_k-p_{k+l}) f(|p_k-p_{k+l}|)+ (p_k-p_{k+N-l}) f(|p_k-p_{k+N-l}|) 
=   4r_0\sin^2{\frac{\theta_l}{2}} f( 2r_0|\sin{\frac{\theta_l}{2}}|) e^{i (\omega_0 t +\theta_k)},
 \end{equation}
 we obtain
 \begin{equation}\label{q_k}
     q_k= \left(i\omega_0 +  F(r_0) \right) r_0 e^{i(\omega_0 t +\theta_k)} ,
 \end{equation}
where $F(r_0)$ is defined in \cref{Fr0}. Substituting \cref{q_k} into \cref{eqhh2} yields
 \begin{equation}\label{eq6-0}
    i\omega_0  \left(i\omega_0  + F(r_0) \right) = \left(M_1- M_2 |i\omega_0 r_0 + r_0F(r_0)|^2 \right) \left(iw_0+ F(r_0) \right)  - k_3  F(r_0). 
 \end{equation}
Equating the imaginary part and the real part of \cref{eq6-0}  gives rise to
\begin{subequations}\label{omega_r}
\begin{align}
        \omega_0^2 &=k_3F(r_0) - F^2(r_0).\\
    M_1&=(1+M_2 k_3 r^2_0) F(r_0).
\end{align}
\end{subequations}
Note that only $M_1$ is determined by the bifurcation parameter. The \cref{omega_r} are solvable when 
\begin{equation}
    M_1 \leq M_{1,c}:= \max_{r_0\geq 0,~F(r_0)\leq k_3} (1+M_2 k_3 r^2_0) F(r_0).
\end{equation}
Therefore, $N$ spots cannot form a rotating ring when $\tau>\frac{1}{k_3}+\frac{M_{1,c}}{k_3^2}$.

\subsection{Stability of rotating N-spot rings}

In this subsection, we  analyze the stability of the ring solution of radius $r_0$ with rotating frequency $w_0$ given by \cref{omega_r}.
We begin by introducing the perturbations to a ring of $N$ spots in the following form
\begin{equation}\label{perbring}
     p_k= r_0 e^{i(\theta_k+\omega_0 t)}(1 + \phi_k), \quad    q_k=  e^{i(\theta_k+\omega_0 t)} \left( i\omega_0 r_0 +r_0F(r_0) +\psi_k \right).
\end{equation}
with $\phi_k$ and $\psi_k$ such that $|\phi_k|,|\psi_k|\ll 1$. Let $l=j-k$, we compute
\begin{equation}
\begin{aligned}
    |q_k|^2 & = ( i\omega_0 r_0 +r_0 F(r_0) +\psi_k )( -i\omega_0 r_0 +r_0 F(r_0)+\bar{\psi}_k )\\
                & = w_0^2r_0^2+r^2_0F^2(r_0) + (r_0F(r_0)-i\omega_0r_0 ) \psi_k +(r_0F(r_0)+i\omega_0r_0) \bar{\psi}_k+|\psi_k|^2.
\end{aligned}
     \end{equation}
Substituting \cref{perbring} into \cref{eqhhode} and neglecting higher order terms leads to the following system:
\begin{subequations}\label{eigenN}
\begin{multline}
    \dot{\phi}_k =-i\omega_0 \phi_k +  \psi_k - \sum_{j\neq k}   (\phi_k- e^{i\theta_l} \phi_j)  f(2r_0|\sin{\frac{\theta_l}{2}}|)\\
   -\sum_{j\neq k}   \frac{ r_0}{4|\sin{\frac{\theta_l}{2}}|} \left[(1-e^{-i\theta_l})(\phi_k+\bar{ \phi}_j)+ (1-e^{i\theta_l})(\bar{\phi}_k+\phi_j)\right] f'(2r_0|\sin{\frac{\theta_l}{2}}|)(1-e^{i\theta_l})
\end{multline}
\begin{multline}
     \dot{\psi}_k =-i\omega_0 \psi_k+  \left( \left[M_1-2M_2(w_0^2r_0^2+r^2_0 F^2(r_0)\right] \psi_k  - M_2\left[ i\omega r_0 +r_0 F(r_0) \right]^2  \bar{\psi}_k \right)\\
    -k_3 \sum_{j\neq k}   (\phi_k- e^{i\theta_l} \phi_j)  f(2r_0|\sin{\frac{\theta_l}{2}}|)
    -\sum_{j\neq k}   \frac{ r_0}{4|\sin{\frac{\theta_l}{2}}|} \left[(1-e^{-i\theta_l})(\phi_k+\bar{ \phi}_j)+ (1-e^{i\theta_l})(\bar{\phi}_k+\phi_j)\right] f'(2r_0|\sin{\frac{\theta_l}{2}}|)(1-e^{i\theta_l}).
\end{multline}
\end{subequations}
With the same notation as in  \cref{sec:3}, we write \cref{eigenN} as
\begin{subequations}\label{eigenNN}
\begin{equation}
    \dot{\phi}_k =-i\omega_0 \phi_k +  \psi_k - \sum_{j\neq k} \left( G_1(\frac{\theta_l}{2}) (\phi_k-e^{i\theta_l} \phi_j)+ G_2(\frac{\theta_l}{2}) ( \bar{\phi}_j-e^{i\theta_l} \bar{\phi}_k  ) \right)
\end{equation}
\begin{multline}
     \dot{\psi}_k =-i\omega_0 \psi_k+  \left( \left[M_1-2M_2(w_0^2r_0^2+r^2_0 F^2(r_0)\right] \psi_k  - M_2\left[ i\omega r_0 +r_0 F(r_0) \right]^2  \bar{\psi}_k \right)\\
   - k_3\sum_{j\neq k} \left( G_1(\frac{\theta_l}{2}) (\phi_k-e^{i\theta_l} \phi_j) + G_2(\frac{\theta_l}{2}) ( \bar{\phi}_j-e^{i\theta_l} \bar{\phi}_k  ) \right)
\end{multline}
\end{subequations}
Assuming that $\phi_k$ satisfies the following relation
\begin{equation}\label{ansatz}
    \phi_k=\xi_+ e^{i m \theta_k}+\xi_- e^{-i m \theta_k},\quad \psi_k=\eta_+ e^{i m \theta_k}+\eta_- e^{-i m \theta_k}, \quad m = 1,\ldots, N.
\end{equation}
Then we can write $\phi_j$ as
\begin{equation}\label{rewrite}
    \phi_j = \xi_+ e^{i m \theta_k} e^{i m \theta_l } +\xi_- e^{-i m \theta_k} e^{-i m \theta_l}, \quad m = 1,\ldots, N.
\end{equation}
Substituting \cref{ansatz} and \cref{rewrite} into \cref{eigenN} and collecting like terms in $e^{im\theta_k},e^{-im\theta_k}$ leads to the system
\begin{subequations}\label{Neigeneq}
\begin{equation}
    i\omega_0 \xi_+ +\dot{\xi}_+ = \eta_+ -  \xi_+ \sum_{j\neq k}  G_1(\frac{\theta_l}{2}) \left(1-e^{i(m+1)\theta_l} \right)- \bar{\xi}_- \sum_{j\neq k} G_2(\frac{\theta_l}{2}) \left(e^{im\theta_l}-e^{i\theta_l} \right) ,
\end{equation}

\begin{equation} \label{xi-}
    i\omega_0 \xi_- +  \dot{\xi}_- = \eta_- - \xi_- \sum_{j\neq k}    G_1(\frac{\theta_l}{2}) \left(1-e^{i(-m+1)\theta_l} \right) -\bar{\xi}_+ \sum_{j\neq k}    G_2(\frac{\theta_l}{2}) \left(e^{-im\theta_l}-e^{i\theta_l} \right) 
\end{equation}

\begin{multline}
     i\omega_0 \eta_+ + \dot{\eta}_+ = \left(  M_1 -2M_2 \left[\omega_0^2 r_0^2 +r_0^2 F^2(r_0)\right] \right) \eta_+  -  M_2\left( i\omega_0 r_0 +r_0 F(r_0) \right)^2 \bar{\eta}_- \\ -k_3\xi_+ \sum_{j\neq k}  G_1(\frac{\theta_l}{2}) \left(1-e^{i(m+1)\theta_l} \right)-  k_3\bar{\xi}_- \sum_{j\neq k} G_2(\frac{\theta_l}{2}) \left(e^{im\theta_l}-e^{i\theta_l} \right) 
\end{multline}

\begin{multline}\label{eta-}
     i\omega_0 \eta_- +  \dot{\eta}_- =   \left(M_1 -2M_2 \left[ \omega_0^2 r_0^2 +r_0^2 F^2(r_0) \right] \right) \eta_- -  M_2\left( i\omega_0 r_0 +r_0F(r_0) \right)^2\bar{\eta}_+   \\   - k_3\xi_- \sum_{j\neq k}    G_1(\frac{\theta_l}{2}) \left(1-e^{i(-m+1)\theta_l} \right) -k_3\bar{\xi}_+ \sum_{j\neq k}    G_2(\frac{\theta_l}{2}) \left(e^{-im\theta_l}-e^{i\theta_l} \right) 
\end{multline}
\end{subequations}
We define
\begin{equation}
    H_1=M_1 -2M_2 \left[ \omega_0^2 r_0^2 + r_0^2 F^2(r_0) \right] ,
\end{equation}

\begin{equation}
    H_2=M_2\left( i\omega_0 r_0 +r_0 F(r_0) \right)^2.
\end{equation}
Using these notations and taking a conjugate of \cref{xi-} and \cref{eta-}  yields
\begin{equation} \label{matrix}
    \begin{pmatrix} \dot{\xi}_+ \\
   \dot{\bar{\xi}}_-\\
   \dot{\eta}_{+}
   \\
   \dot{\bar{\eta}}_{-}
   \end{pmatrix}=
   \begin{pmatrix}
 -I_1(m) -i\omega_0 &  -I_2(m) & 1 & 0 \\
   -I_2(m)& -I_1(-m) +i\omega_0 & 0&1\\
  -k_3 I_1(m)  &  -k_3 I_2(m)&H_1-i\omega_0&-H_2\\
   -k_3 I_2(m)& -k_3 I_1(-m) &-\bar{H}_2&H_1+i\omega_0
    \end{pmatrix}
      \begin{pmatrix} \xi_+ \\
   \bar{\xi}_-\\
   \eta_{+}
   \\
   \bar{\eta}_{-}
   \end{pmatrix}
\end{equation}
Let  $\begin{pmatrix} \xi_+ \\
   \bar{\xi}_-\\
   \eta_{+}
   \\
 \bar{\eta}_{-}
   \end{pmatrix}=e^{\lambda t}
   \begin{pmatrix}
   a_1\\
   a_2\\
   a_3\\
   a_4
   \end{pmatrix}$, then $\lambda$ is an eigenvalue of the $4\times 4$ matrix in \cref{matrix}. The eigenvalue $\lambda$ must satisfy:
   \begin{equation}
       \lambda a = M a
   \end{equation}
 where 
 \begin{equation}
     M= 
   \begin{pmatrix}
 -I_1(m) -i\omega_0 &  -I_2(m) & 1 & 0 \\
   -I_2(m)& -I_1(-m) +i\omega_0 & 0&1\\
  -k_3 I_1(m)  &  -k_3 I_2(m)&H_1-i\omega_0&-H_2\\
   -k_3 I_2(m)& -k_3 I_1(-m) &-\bar{H}_2&H_1+i\omega_0
    \end{pmatrix}
 \end{equation}
 In this way, the study of stability of an rotating $N$-spot ring solution decouples into the study of individual Fourier modes. We conclude the results in \cref{prop3}.

 \section{Numerical validation} \label{sec:5}
 
Solutions to the PDE system~\cref{rd} are computed numerically with Fourier-type pseudo-spectral method in space and the MATLAB subroutine ``ode113'' for the time evolution. The MATLAB codes are available at \url{https://github.com/KaleonXie/Interaction-of-Spots-with-Oscilltory-Tails}. To examine the formation of an $N$-spot ring, we pick parameters that provide spot solutions with tails that decay to the homogeneous background state. This results in an interaction function which also shows oscillatory behaviour with attractive $(f(d) > 0)$ and repulsive $(f(d) < 0)$ regions of interaction. There are theoretically an endless number of $N$-spot rings. However, when the radius of the ring exceeds the second binding radius, the interactions between spots are so weak that they are undetectable in numerical calculations. As a result, we concentrate on the numerical analysis of the $N$-ring with the first and the second smallest binding radius.

The numerical verification of the stability requires a significant amount of computational work, since we must ensure that any interaction between neighboring spots has completely vanished.  We note that long-time simulations are necessary to check the stability of the stationary, traveling, and rotating $N$-spot rings. In order to capture the slow dynamics, we let the system evolve until $t=40000$. Numerous numerical experiments are conducted to verify the results. In all of the numerical computations below, we choose the parameters:
\begin{equation}
D_u=1.1\times 10^{-4},~D_w=9.64\times 10^{-4},~k_1=1.01,~k_3=0.3,~k_4=1,~\kappa=-0.1.
\end{equation}
The spatial discretization is $256\times 256$ in a square $[-1,1]\times[-1,1]$.

\begin{figure}[!htb]
\centering
    \begin{subfigure}[b]{0.48\textwidth}
       \centering
         \includegraphics[width=\linewidth]{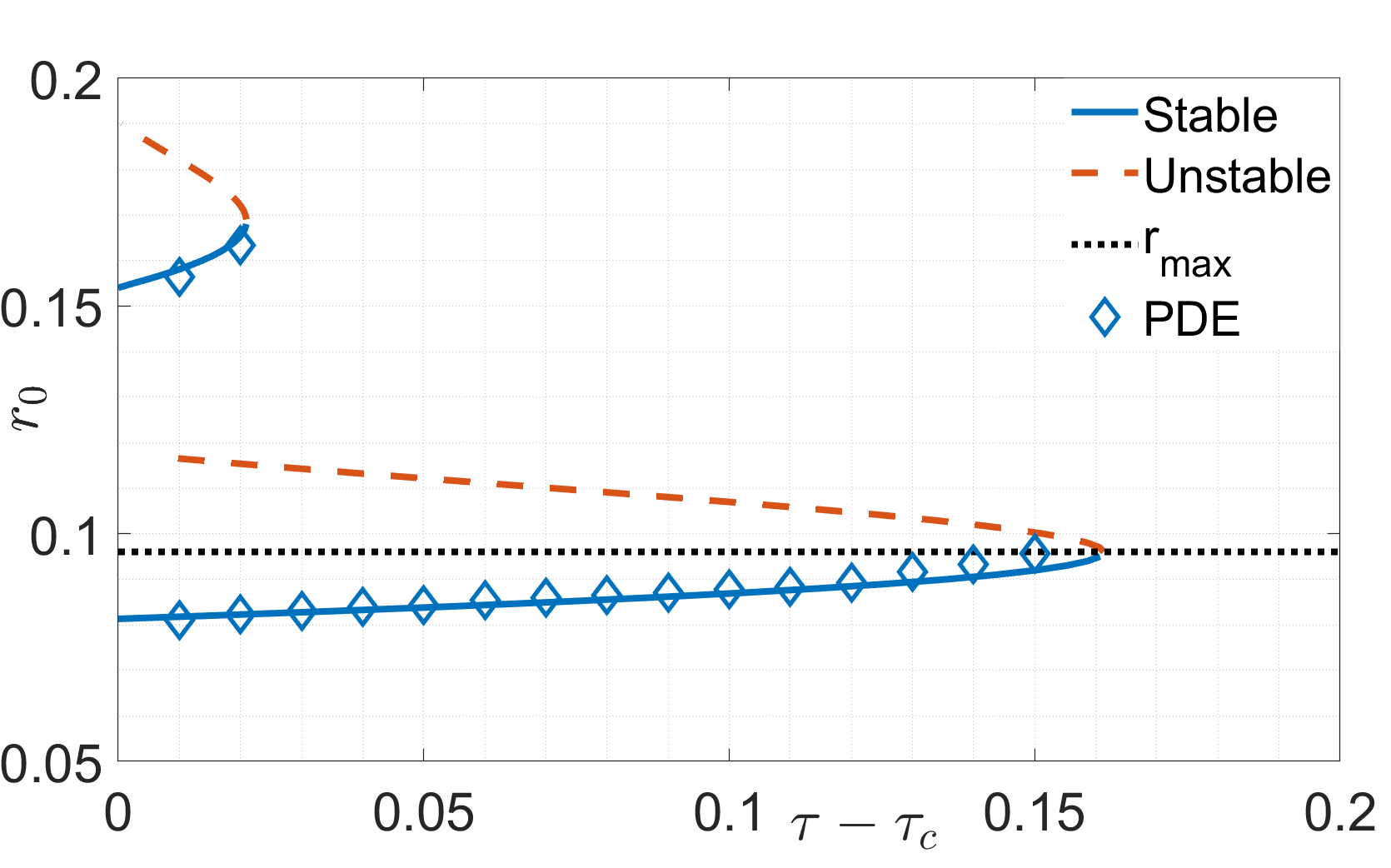}
         \caption{Rotating two-spot ring.}\label{fig:2-spot}
    \end{subfigure}
    \begin{subfigure}[b]{0.48\textwidth}
        \centering
         \includegraphics[width=\linewidth]{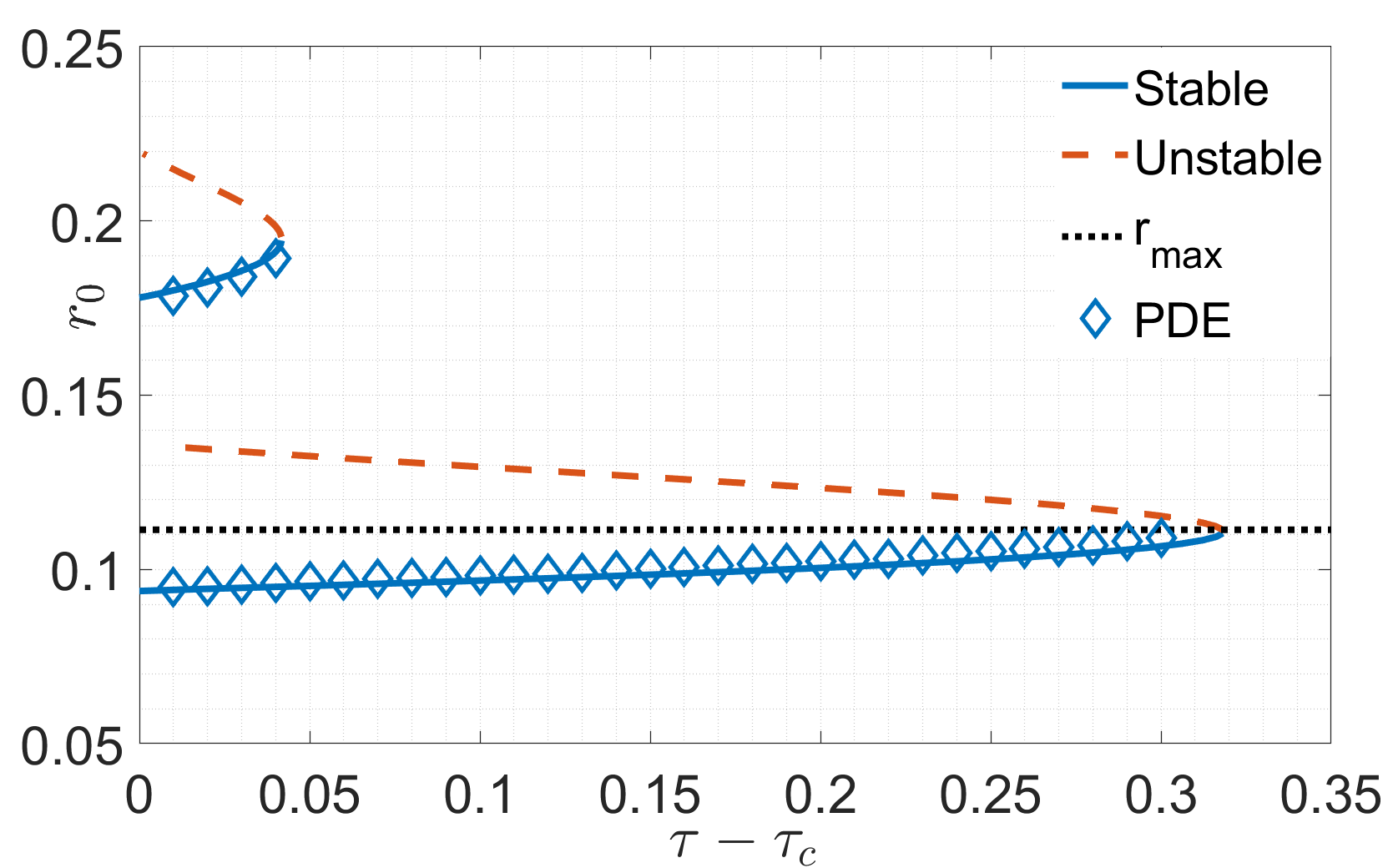}
         \caption{Rotating three-spot ring.}\label{fig:3-spot}
    \end{subfigure}   
    \begin{subfigure}[b]{0.48\textwidth}
       \centering
         \includegraphics[width=\textwidth]{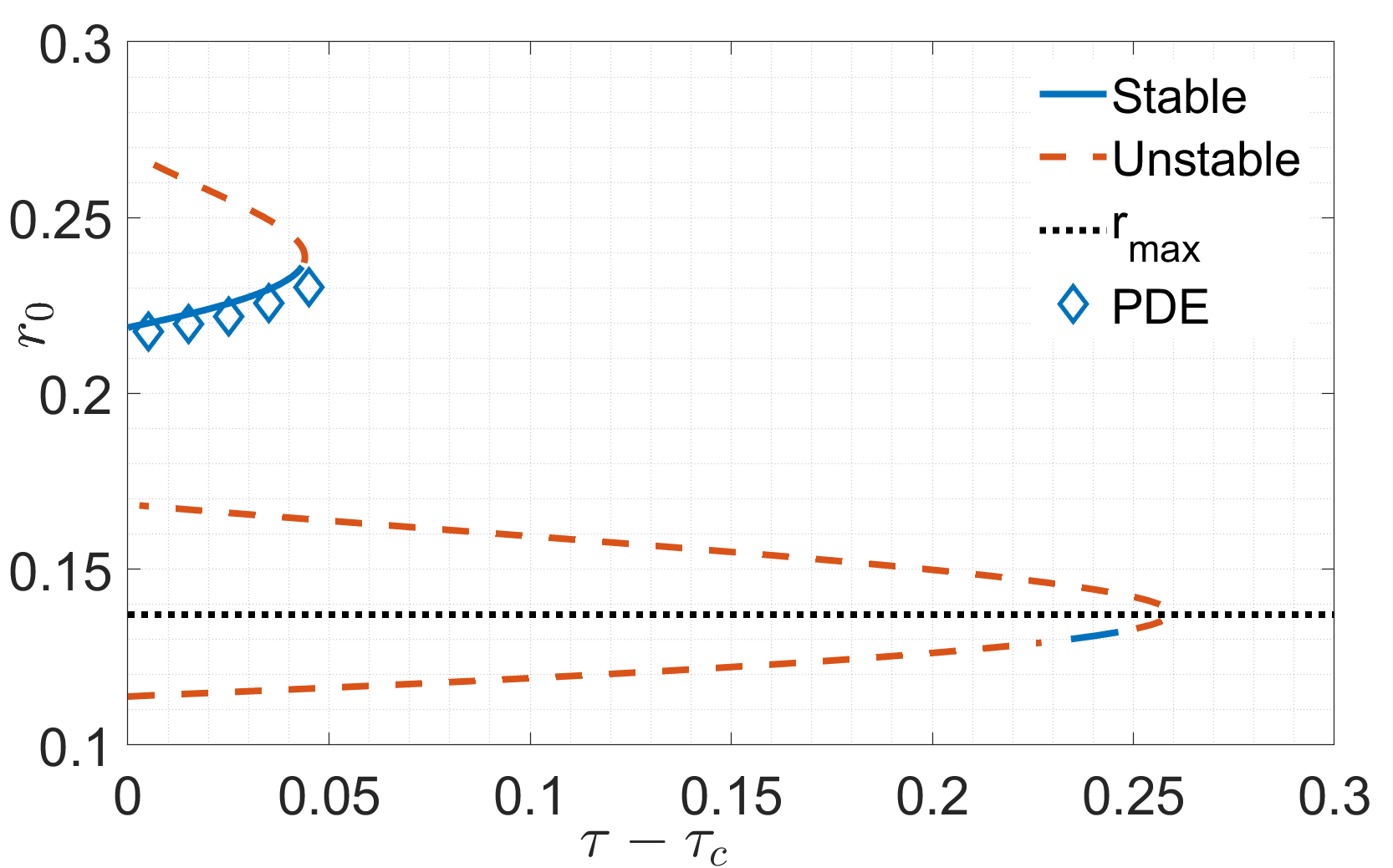}
         \caption{Rotating four-spot ring.}\label{fig:4-spot}
    \end{subfigure}     
    \begin{subfigure}[b]{0.48\textwidth}
    \centering
         \includegraphics[width=\textwidth]{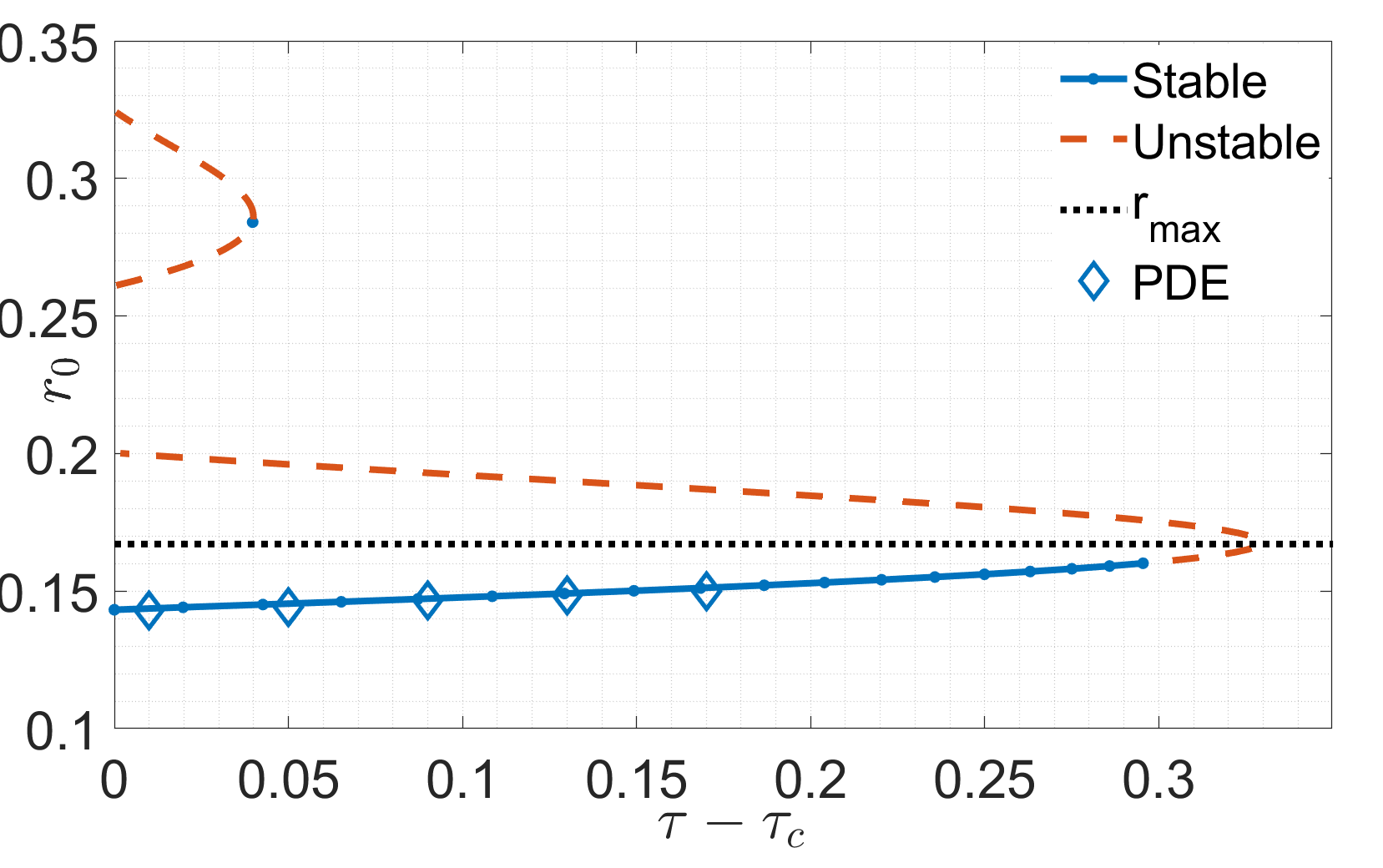}
         \caption{Rotating five-spot ring.}\label{fig:5-spot}
    \end{subfigure}
    \caption{ Radius of rotational $N$-spot ring for $N=2,3,4,5$ as a function $\tau-\tau_c$.  Other parameters are given in the caption of \cref{N-spot}. Solid (dash) lines are stable (unstable) states predicted from \ref{prop3}. Diamonds indicate the mean radii from the simulations of PDE. Dot horizontal lines $r_0=r_{max}$ correspond to the radius of maximal possible rotating state.}
    \label{fig:2-5-spot}
\end{figure}

\begin{Experiment}
Stationary $N$-spot rings below the threshold.
\end{Experiment} 

Let $\tau=0.1<\tau_c$. We first place the initial spots on a ring with the first or second binding radius according to the root of \cref{r0}. Then,  to obtain a stationary $N$-spot ring, we run the PDE simulations until the $L_\infty$ norm of the difference between the state at two successive times, $T$ and $T+10$,  is smaller than some tolerance ($10^{-7}$ in our simulation).  After obtaining a stable stationary $N$-spot ring, we add a minor perturbation to the stationary state and let the system evolve until the difference between two subsequent states is smaller than the tolerance again.  The stability predicted by our criterion is in good agreement with the ODE and PDE simulations. We note that $N$-spot rings with the first binding radius are unstable when $N=4,7$, whereas $N$-spot rings with the second binding radius are always stable.  \cref{table:ss-N-s} summarizes the stability results for stationary $N$-spot rings.
   
  There is a difficulty that can arise in obtaining an $N$-spot ring with the first binding radius. If the oscillatory tail has a large amplitude, a superposition of $N$ spots at the ring center can be sufficient to ignite additional spots.  For the parameters we use, numerical simulations show that a new spot emerges at the center when $N=6$. We can obtain an $N$-spot ring for $N\neq 6$, the stability of which can be predicted by the reduced ODE \cref{ordinary}.  For ensembles of $N$-spot rings with the second binding radius, the superposition of the tail at the center falls behind the decay of the tail, allowing for the development of an $N$-spot ring. \cref{eqhhode} can be used to more precisely explain the dynamics.
    
\begin{table}[!htb]
\begin{center}
\begin{tabular}{|c|c|c|c|c|c|c|c|}
\hline
 $N$ & 2 & 3 & 4 & 5 & 6  & 7 & 8  \\  \hline
 ODE BS I & stable & stable & unstable & stable & stable & unstable & stable \\  \hline
PDE BS I  &stable  &stable  &unstable  & stable &N.A.  &unstable &stable \\  \hline
ODE BS II & stable & stable & stable & stable & stable  & stable & stable  \\  \hline
PDE BS II& stable & stable & stable & stable &   stable  & stable & stable  \\  \hline
\end{tabular}
\end{center}
\caption{Stability of the stationary N-spot rings for $N=2,\ldots, 8$ when $\tau=0.1<\tau_c$. N.A. means not available here. ODE(PDE) BS I and BS II are referred to as the $N$-spot ring state with the first and second binding radii in ODE(PDE) simulation. Movies of simulations are provided in the supplement material.}
\label{table:ss-N-s}
\end{table}    
\begin{Experiment}
Traveling $N$-spot rings near the threshold.
\end{Experiment} 

 To verify the existence and stability of traveling $N$-spot rings, we start with the stationary $N$-spot ring obtained at $\tau=0.1$ and add a small initial uniform velocity to it. Then we let the simulation run until $t=40000$. \cref{table:t-ss-N-s} summarizes the stability results for traveling N-spot rings at $\tau=\tau+0.01$.   We remark that the stable traveling spot obtained in PDE simulation in \cref{table:t-ss-N-s} is referred to as a traveling spot with a fixed speed without taking the direction into account. Namely, a two-spot ring traveling horizontally is equivalent to a two-spot ring moving vertically. In the PDE simulation, a two-spot ring will move eventually along its longitude direction, which cannot be predicted by the reduced ODE. Higher order terms are necessary to explain this transition.
 
     When $N>7$, even though the ODE simulation shows stable traveling rings, the PDE simulation may give different results. In the PDE simulation, once a $N$-spot ring starts to travel, there occurs a symmetry-breaking from the regular ring shape, which may initiate the higher-order interaction. The deformation usually consists of a small elongation along the traveling direction (or equivalently, a tiny shrinkage in the orthogonal direction), which makes the distance slightly shorter between spots located alongside. Then higher-order interaction (of the attractive type) is evoked.

  \begin{remark}
     As our second order reduced ODE system \cref{eqhh} differs from the second order ODE model in swarming \cite{albi2014stability} by one term, it is interesting to see the difference between them. In the first and second order ODE models studied in swarming, the stability of stationary $N$-spot ring in the first and second order models are equivalent. This is no longer true in our reduced system. A six-spot ring with the second binding radius is stable when it is stationary but unstable when it starts to move. 
  \end{remark}

    \begin{table}[!htb]
\begin{center}
\begin{tabular}{|c|c|c|c|c|c|c|c|}
\hline
 $N$ & 2 & 3 & 4 & 5 & 6  & 7 & 8  \\  \hline
 ODE BS I & stable & stable & unstable & stable & stable & unstable & stable \\  \hline
PDE BS I  &stable  &stable  &unstable  & stable &N.A.  &unstable &unstable \\  \hline
ODE BS II & stable & stable & stable & unstable & unstable  & stable & stable  \\  \hline
PDE BS II& stable & stable & stable & unstable &   unstable & unstable & unstable \\  \hline
\end{tabular}
\end{center}
\caption{Stability of the traveling N-spot rings for $N=2,\ldots, 8$ when $\tau=\tau_c+0.01$. Notations are the same as \cref{table:ss-N-s}. Movies of simulations are provided in the supplement material.  }
\label{table:t-ss-N-s}
\end{table}    
    
\begin{Experiment}
Rotating $N$-spot rings near the threshold.
\end{Experiment} 

To verify the existence and stability of rotating $N$-spot rings, we start with the stationary $N$-spot ring obtained at $\tau=0.1$ and add a small initial rotational velocity to it. Then we let the simulation run until $t=40000$ to obtain a possible rotating state. In the middle of simulation, $t=20000$, we also add a small velocity to one spot on the ring to check the stability of this rotating state.  A detailed investigation of rotational $N$-spot rings for $N=3,4,5$ is depicted by \cref{fig:2-5-spot}.  \cref{table:r-ss-N-s} summarizes the stability results for rotating N-spot rings at $\tau=\tau_c+0.01$.  Rotating $N$-spot rings with the second binding radius are stable for $N=2,3,4$ and unstable for $5 \leq N\leq 8$, in agreement with the prediction from \cref{prop3}.
    
     \begin{table}[!htb]
\begin{center}
\begin{tabular}{|c|c|c|c|c|c|c|c|}
\hline
 $N$ & 2 & 3 & 4 & 5 & 6  & 7 & 8  \\  \hline
 ODE BS I & stable & stable & unstable & stable & stable & unstable & stable \\  \hline
PDE BS I  &stable  &stable  &unstable  & stable &N.A.  &unstable & unstable \\  \hline
ODE BS II & stable & stable & stable & unstable & unstable  & unstable & unstable  \\  \hline
PDE BS II& stable & stable & stable & unstable &   unstable  &unstable & unstable  \\  \hline
\end{tabular}
\end{center}
\caption{Stability of the rotating N-spot rings for $N=2,\ldots, 8$ when $\tau=\tau_c+0.01$. Notations are the same as \cref{table:ss-N-s}.  Movies of simulations are provided in the supplement material. }
\label{table:r-ss-N-s}
\end{table}

   Stability examination of the rotating four-spot ring  in \cref{fig:4-spot} reveals that a rotating ring can be stable despite its stationary counterpart being unstable.   A stationary four-spot ring with the first binding radius under a small perturbation will contract and eventually take the form of a rhombus. On the other hand, a rotating ring with a high rotational speed will cause the ring's radius to increase to the point where the centripetal force produced by other spots cannot maintain the revolution. There will be an intermediate regime where the shrinkage balances the expansion, resulting in a stable rotating ring. However, it is only visible in the ODE simulations. In PDE simulations, the rotating four-spot ring with the first binding radius is always unstable.

    Ignition of a new spot in the center of the five-spot ring can be observed as we increase the parameter $\tau$.  As $\tau$ grows, the radius of the ring expands, and the superposition in the center climbs to levels exceeding the igniting threshold. Thus, a rotating five-spot ring is unstable in PDE simulations at speeds significantly less than those anticipated by ODE. see \cref{fig:5-spot}.

\section{Conclusions and Outlook  } \label{sec:6}
Our basic question in the pattern formation field is that how a group of spots with oscillatory tails interact and what kind of ordered state they form asymptotically. Are they different from those with monotone tails? For the monotone case, it is either repulsive or attractive. However, there are infinitely many possibilities for the oscillatory case, but practically only the first few interactions matter due to exponentially decaying. We are interested in a spontaneous formation of ring patterns, which serve as fundamental building blocks for complex dynamics and are not driven by boundary conditions or the shape of the domain \cite{kolokolnikov2022ring}. In this article, we have investigated the stationary and moving ring solutions of a RD system both analytically and numerically. When the reaction rate is below the threshold, the slow dynamics due to the spot-spot interaction can be described by a first-order ODE system, through which we are able to determine the existence and stability of a stationary $N$-spot ring.  When the reaction rate is slightly above the threshold, self-propelled motions of the spots are induced by the drift instability, leading to the traveling and rotational motions of $N$-spot rings. The dynamics of the PDE system can be described by a set of ODEs that have the same bifurcating structure as the original system.  We then demonstrate the existence of traveling and rotating $N$-ring solutions for the reduced ODE system,  the stability of which is determined by the eigenvalues of $N$ 
matrices of $4\times 4 $ size. Our analytical results are validated by numerical simulations of PDE and ODE systems.

When we start from a general initial arrangement of many spots in the PDE simulation, we observe many voids (homogeneous regions surrounded by spots) in each transient cluster of spot after the initial transient. Voids appear during the interacting process and some of them remain, but some of them can be relaxed through generating new spots. The dynamics of ring patterns clearly illustrate these processes. Instability of $N=4,6,7$ for stationary $N$-spot rings with the first binding radius is a key to understanding these dynamics, we refer the readers to the relevant movies in the supplement materials. In the four-spot ring's simulation, two spots in the diagonal migrate slowly toward the center and occupy the void. In contrast, in the simulation of seven-spot ring, the ring gradually deforms but the void in the center persists. Both of these two dynamics are attributed to higher-order-term interactions originated from non-neighbouring spots. The void in the center is eased in the simulation of six-spot ring.  The emergence of a spot from a void is due to the contamination of the activator (a significant disruption of the homogeneous state through tail overlapping).  The above dynamics are spontaneous and we regard them as ``self-repairing” forces to change one state to a more stable one. Because these instabilities do not occur in spots clustered in a convex shape with the nearest binding distance, we speculate that \textit{“any convex shape of cluster in which any two spots are bound with the first binding distance is stable”}.

We have illustrated the interaction of spots under a special parameter setting for PDE: $D_v =0$ and $\theta=0$, which have been used in a series papers to study the interaction of dissipative solitons, see \cite{liehr2013dissipative}. The advantage of this choice is that the PDE is simple enough to produce traveling solitons, and deriving the reduced system near the drift bifurcation is relatively straightforward.  We note that it is feasible to derive similar form of reduced models for $D_v \neq 0$ and $\theta \neq 0$ by following the approaches in \cite{ei2002motion}. Thus the analysis can be easily extended to general three-component systems.

There are many open problems left to be explored for the localized spots with oscillatory tails. Partial list is as follows:
\begin{itemize}

    \item  \textbf{Compact arrangement of multiple spots}: Since spots can link to one another with different binding distances, numerous stationary, stable states may be established for spots with oscillatory tails.   One intriguing pattern is the dense arrangement of spots, where spots form a compact structure with the shortest binding distance and no voids appear inside. Special instances of these are the two-spot and three-spot rings with the smallest binding radius.  It is interesting is to investigate the number of stable stationary compact configurations for a fixed number of spots. \cref{5-spotv} depicts a variety of  stable $5$-spot contours. Similar cluster patterns have also been observed in plane gas-discharge experiments \cite{nasuno2003dancing}. 
    
    \begin{figure}[!htb]
         \centering
         \includegraphics[width=\textwidth]{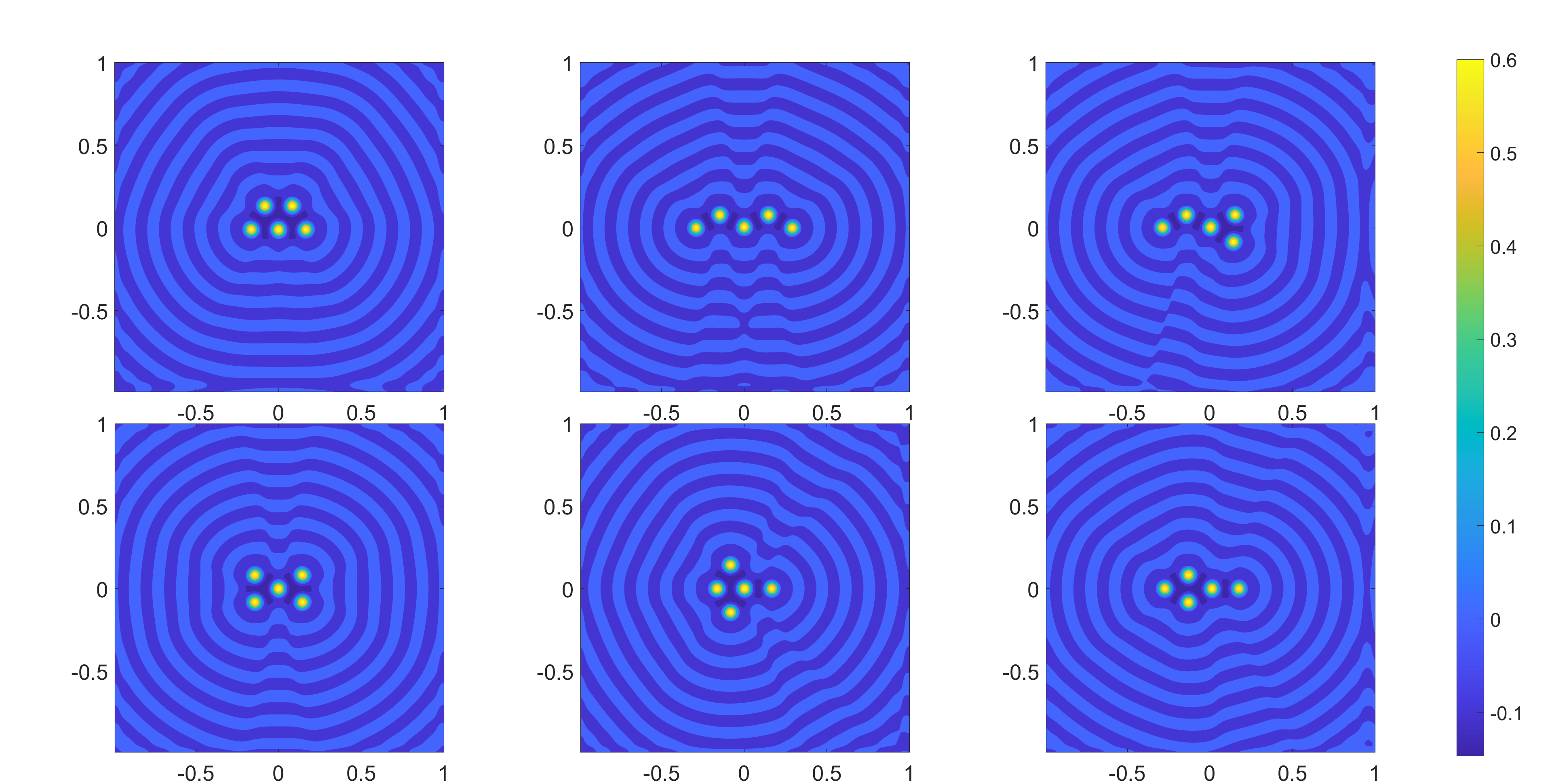}
         \caption{ Various stable five-spot contours of $u-u_c$ at $\tau=0.1$. Other parameters are the same as in \cref{N-spot}. }
         \label{5-spotv}
\end{figure}
    
    \item  \textbf{Direction of moving spots:}  Numerical simulations show that all spots may march in a preferred direction at the final stage, despite their beginning orientations being different.  This cannot be explained by the reduced ODE system and warrant further investigation.

\item  \textbf{Collision:} It is very interesting to explore the colliding dynamics of traveling spots with oscillatory tails. Numerical simulations have revealed new dynamic behaviors, not seen in the simulation of spots with monotone tails. Except for fusion and annihilation, we also observe the creation of new spots when multiple spots move toward one point. A typical scenario is the emergence of a new spot in the center of a six-spot ring with the first binding radius. It is due to the existence of the ``scattors",  see \cite{nishiura2005scattering,nishiura2003dynamic}. As the homogeneous state is locally stable,  a perturbation with an amplitude above some threshold is necessary to get a locally radially symmetric spot solution. This observation suggests that there is a smaller spot of saddle type known as a scattor that plays an important role in determining the evolution's final state.  Another new phenomenon is the formation of rotating ring. Two or three traveling spots colliding with off-center may form a binary (triply bonded) star. The final rotating ring state have been discussed in this paper. However, it is not clear under what conditions the rotating ring can develop.

\item  \textbf{1D pulses with oscillatory tails:} There are many types of stationary bound states constituted by pulses with oscillatory tails in 1D. In contrast with 2D spots, a 1D pulse only interacts with its neighbouring pulses. Thus, an $N$-pulse bound state consists of pulses with different binding distances to their neighbours.  The existence and stability of various bound states remains to be systematically studied.

\item \textbf{Interplay between different modes:} The transition of a stationary single spot to a rotating spot has been reported and studied in  \cite{teramoto2009rotational,xie2017moving}. The rotational motion in \cite{teramoto2009rotational} is boundary-free and caused by the interplay of the translational and splitting modes. While the rotating spot in \cite{xie2017moving} requires a Neumann boundary condition on a disk domain and is triggered by the Hopf bifurcation associated with the translational mode.  The instability of these two distinct types of revolving spot remains to be explored. 

\item \textbf{Ring Structure in other systems:} Recently, a ring of spikes for the Schnakenberg model inside either a unit disk or an annulus has been considered in \cite{kolokolnikov2022ring}.  They have shown that a ring of eight or less spikes is stable inside a disk. However, for Schnakenberg model, the ring of spikes only exists in certain special domains, since the spot-spot interaction is controlled by another component that does not localize and is much stronger and domain-dependent. It is intriguing to investigate the possible paths of these spikes when a $N$-spike ring becomes unstable with respect to the drift mode.

\end{itemize}

\section*{Acknowledgments} 
S.X. and Y.N. acknowledge partial support by the Council for Science, Technology and Innovation (CSTI), Japan,
Cross-Ministerial Strategic Innovation Promotion Program (SIP), Japan, ‘Materials Integration’ for Revolution-
ary Design System of Structural Materials. Y.N. gratefully acknowledges the support by JSPS KAKENHI Grant number JP20K20341. Y.N. also thanks Professor Kei-Ichi Ueda for valuable comments on the reduced ODE system.

\bibliographystyle{unsrt}

\bibliography{ref}
\end{document}